\documentclass[11pt,fleqn]{article}
\RequirePackage[OT1]{fontenc}
\RequirePackage{amsthm,amsmath}
\RequirePackage[round]{natbib}
\RequirePackage[colorlinks,citecolor=blue,urlcolor=blue]{hyperref}
\usepackage{pgf,tikz}\usepackage{mathrsfs}\usetikzlibrary{arrows}
\usepackage{natbib}
\usepackage{fullpage}
\usepackage{amsmath}
\usepackage{amssymb}
\usepackage{graphicx}
\usepackage{dsfont}
\makeatletter
\def\captionof#1#2{{\def\@captype{#1}#2}}
\makeatother
\def\1{\mbox{\bf 1}}
\def\R{\mathbb{R}}

\def\N{\mathbb{N}}
\def\P{\mathbb{P}}
\def\E{\mathbb{E}}
\def\L{\mathbb{L}}

\def\R{\mathbb{R}}

\def\Z{\mathbb{Z}}

\newtheorem{theo}{Theorem}
\newtheorem{lem}{Lemma}
\newtheorem{prop}{Proposition}

\newtheorem{cor}{Corollary}
\newtheorem{Def/Prop}{Definition-Proposition}

\newcounter{exos}
\renewcommand\theexos{\arabic{exos}}

\newcounter{prob}
\renewcommand\theprob{\arabic{prob}}

\begin{document}
\author{Lionel Truquet \footnote{UMR 9194 CNRS CREST, ENSAI, Campus de Ker-Lann, rue Blaise Pascal, BP 37203, 35172 Bruz cedex, France. {\it Email: lionel.truquet@ensai.fr}.}
 }
\title{Coupling and perturbation techniques for categorical time series}
\date{}
\maketitle

\begin{abstract}
\noindent
We present a general approach for studying 
autoregressive categorical time series models with dependence of infinite order and defined conditional on an exogenous covariate process. 
To this end, we adapt a coupling approach, developed in the literature for bounding the relaxation speed of a chain with complete connection and from which we derive a perturbation result for non-homogenous versions of such chains.  
We then study stationarity, ergodicity and dependence properties of some chains with complete connections and exogenous covariates. As a consequence, we obtain a general framework for studying some observation-driven time series models used both in statistics and econometrics but without theoretical support.
\end{abstract}
\vspace*{1.0cm}

\footnoterule

\section{Introduction}
Categorical time series are widely encountered in various fields. For instance, in climate analysis, \citet{Guanche} studied the dynamic of weather types, \citet{Hao} the prediction of drought periods. In finance, \citet{russell} or \citet{rydberg} studied the dynamic of price movements. In economics, \citet{kauppi} consider the prediction of recession periods. 
Several type of models used for modeling categorical time series can be found in the survey of \citet{Fokianos}. 
Though lots of time series models have been developed in the literature, it is difficult to find a general framework for which inclusion of exogenous covariates is mathematically justified. This is one of the important differences between the theoretical results found in time series analysis and the models used by the practitioners which most of the time, are based on exogenous covariates. 
A notable exception is the contribution of \citet{Kaufmann} who considered estimation in autoregressive logistic
type models when deterministic regressors are included in the dynamic. More recently, \citet{FT1} considered general Markov models with random covariates.
However, most of the categorical time series models used in practice are "observation-driven" (see below for a definition), especially in econometrics. \citet{FT1} also considered this class of non-Markovian processes but without covariates and it seems that a general approach for studying a wide class of categorical time series models with exogenous covariates is still not available. 
In this paper, we provide such a framework by using a formalism introduced for studying a general class of finite-state stochastic processes, the chains with complete connections. 
These processes, initially considered by \citet{DF}, have an interest in probability theory, statistical mechanics or ergodic theory.  
See in particular \citet{Harris}, \citet{Ios}, \citet{Bres}, \citet{Bres2}, \citet{GF} and \citet{comets} for many of their theoretical properties. Chains with complete connections also contain stochastic chains with memory of variable length as a special case, the latter 
class, initially introduced by \citet{Ris} for data compression, has also applications in linguistic, see \citet{Galves} or for protein classification, see for instance \citet{prot}.

In this paper, we also consider such chains with complete connections but defined conditional on a covariate process. More precisely, we want to study stochastic processes $(Y_t)_{t\in\Z}$ defined by 
\begin{equation}\label{prob}
\P\left(Y_t=y\vert Y_{t-1}^{-},X_t^{-}\right)=q\left(y\vert Y_{t-1}^{-},X_t^{-}\right),\quad y\in E,
\end{equation}  
where $(X_t)_{t\in\Z}$ is a covariate process taking values in $\R^d$, $E$ is a finite set and $q$ is a transition kernel. We will extensively use the notation $x_t^{-}=(x_t,x_{t-1},\ldots)$ for a sequence $(x_t)_{t\in\Z}$.
 Without additional assumptions on the two processes $X$ and $Y$, (\ref{prob}) is difficult to study theoretically. We will assume further that 
\begin{equation}\label{prob2}
\P\left(Y_t=y\vert Y_{t-1}^{-},X_t^{-}\right)=\P\left(Y_t=y\vert Y_{t-1}^{-},X\right),\quad X:=(X_t)_{t\in\Z}.
\end{equation}
If condition (\ref{prob2}) is satisfied, $(Y_t)_{t\in\Z}$ is, conditional on $X$, a time-inhomogenous chain with complete connections and transition kernels $\left\{q\left(\cdot\vert\cdot,X_t^{-}\right): t\in\Z\right\}$. 
Condition (\ref{prob2}) also means that $Y_t$ is independent of $(X_{t+1},X_{t+2},\ldots)$conditional on $\left((Y_{j-1},X_j)\right)_{j\leq t}$. 
In econometrics, the latter conditional independence assumptions is called strict exogeneity. Initially introduced by \citet{Sims} for linear models, the concept of strict exogeneity was extended by \citet{Chamb} to categorical time series. 
\citet{Chamb} also showed that under additional regularity conditions, 
strict exogenity is equivalent to non Granger causality, which means that $X_{t+1}$ is independent of $Y_t,Y_{t-1},\ldots,$ conditional on $X_t,X_{t-1},\ldots$. This roughly means that the covariate process evolves in a totally autonomous way and that, given all the information available up to time $t$, past values of the outcome will not influence future values of the covariates.
Let us also mention that such strict exogeneity condition is a standard assumption in Markov-switching models, for which the dynamic of the time series under study is defined conditional on an unobserved Markov chain $X$.     
In probability theory, this exogeneity notion appears implicitly in the literature of stochastic processes in random environments. 
Finite-state Markov chains in random environments are a particular case of stochastic processes satisfying (\ref{prob}) and (\ref{prob2}). They are studied for instance in \citet{Cog} and \citet{Kif} but no result seems to be available for chains with complete connections.
Strict exogeneity has of course some limitations for time series analysis, it is a rather strong assumption. However, it is easier to formulate a general theory in this context, other conditional distributions such as $X_{t+1}\vert Y_t^{-},X_t^{-}$ need not to be specified. 

Stochastic processes defined by (\ref{prob}) are of course of theoretical interest but 
for applications to time series analysis, one of the challenging problem is to find parsimonious 
versions of (\ref{prob}). One important class of models are called observation-driven, following the classification proposed by \cite{Cox}. For model (\ref{prob}), an observation-driven model is obtained assuming that $q\left(\cdot \vert Y_{t-1},X_t^{-}\right)=q\left(\cdot\vert \mu_t\right)$ with 
\begin{equation}\label{od}
\mu_t=G\left(\mu_{t-1},\ldots,\mu_{t-q},Y_{t-1},\ldots,Y_{t-p},X_t\right).
\end{equation}
Without exogenous covariates, observation-driven models were widely studied, in particular for count time series. See in particular \citet{Fok},\citet{Neu},\citet{EJS},\citet{DDM}.
This models are mainly studied using Markov chain techniques due to the Markov properties of the process $(Y_t,\mu_t)$. However, as pointed out in \citet{EJS} or \citet{DDM}, for discrete time series, such Markov chains do not satisfy irreducibility properties. In particular the latent variable $\mu_t$ is not discrete and not necessarily absolutely continuous. More sophisticated techniques have then been developed to study existence of stationary distributions. Such contributions are often limited to the case $p=q=1$ and do not consider the problem of exogenous covariates. In contrast, for the special case of categorical time series, 
one can develop a much more general approach, considering observation-driven models as a particular case of infinite dependence. This approach was recently used by \cite{FT1}.
However, inclusion of exogenous covariates is a more tricky problem and has not been considered before for model (\ref{od}) or (\ref{prob}). More generally, despite its fundamental importance for practical applications, the problem of covariates inclusion is often ignored in the time series literature, except for linear models. In Section \ref{4}, we make a review of many observation-driven models proposed in econometrics for the study of categorical time series and that can be studied under our general framework.

A crucial point for studying our models is to control how fast the process $(Y_t)$ in (\ref{prob}) loses memory of its initial values. For homogenous chains, \citet{Bres} developed a nice result based on the maximal coupling. We will adapt their result to our context, which will be crucial for defining our models and studying many of their properties.  

Another important problem addressed in this paper concerns dependence properties of the process, which are essential to control the behavior of partial sums. While chains with complete connections satisfies $\phi-$mixing properties under rather general assumptions (see \citet{FT1}), finding dependence properties for the joint process $(Y_t,X_t)$ in (\ref{prob}) is quite challenging. For the example of observation-driven models, such properties are crucial to control the behavior of partial sums of type $\frac{1}{n}\sum_{t=1}^n f\left(Y_t^{-},\mu_t^{-}\right)$. In this paper, we explain how to get $\beta-$mixing properties and $\tau-$dependence (see Section \ref{5} for a definition) for this joint process. 
To this end, we will use a coupling approach. If a "good" coupling for the covariate process $X$ exists, one can define a coupling of $Y$ conditional on $X$, with two paths having different transition kernels that will be adjacent at infinity. This is why we will derive in Section \ref{2} a perturbation result for chains with complete connections, obtained via coupling. Such a result also has an independent interest.

The paper is organized as follows. In Section \ref{2}, we state a general result for non-homogenous chain with complete connections. In particular, we generalize a result of \citet{Bres} for controlling the relaxation speed of such chains and we also compare the dynamic of two such chains possessing different transition kernels. 
In Section \ref{3}, we give some conditions on the transition kernel $q\left(\cdot\vert\cdot\right)$ that guaranty existence and uniqueness of a stationary and ergodic solution for the problem (\ref{prob}). Many examples are given in Section \ref{4}, with a detailed treatment of some observation-driven models used in the econometric literature.
Section \ref{5} is devoted to the dependence properties of the solution, absolute regularity or $\tau-$dependence. We mention several possible applications of our results in statistics in Section \ref{6}. Finally, several auxiliary lemmas for the proofs of our results are collected in an Appendix.

\section{Perturbation of chains with complete connection}\label{2}
We denote by $\N$ the set of natural integers $\{0,1,\ldots\}$ and $\N^{*}=\N\setminus\{0\}$.For a finite set $F$, we will denote by $\mathcal{P}(F)$ is the set of all subsets of $F$.
Moreover if $\nu_1$ and $\nu_2$ are two probability measures on $F$, the total variation distance between $\nu_1$ and $\nu_2$ is defined by 
$$d_{TV}\left(\nu_1,\nu_2\right)=\frac{1}{2}\sum_{f\in F}\left\vert \nu_1(f)-\nu_2(f)\right\vert.$$
We remind that we have the following dual expression 
$$d_{TV}\left(\nu_1,\nu_2\right)=\inf\left\{\P\left(U\neq V\right): U\sim \nu_1, V\sim \nu_2\right\}.$$
For $y,\overline{y}\in E^{\N}$ and a positive integer $m$, we  write $y\stackrel{m}= \overline{y}$ if $y_i=\overline{y}_i$ for $0\leq i\leq m-1$.

\subsection{A general result}
Throughout the section, we will denote by $E$ a finite set. Let $\left(\mathcal{X},\mathcal{B}\left(\mathcal{X}\right)\right)$ be a Polish space.
For any $x\in \mathcal{X}$, we consider two sequences $\left(q^{x}_t\right)_{t\in\Z}$ and $\left(\overline{q}^{x}_t\right)_{t\in\Z}$ of probability kernels from $\left(E^{\N},\mathcal{P}(E)^{\otimes\N}\right)$ to $\left(E,\mathcal{P}(E)\right)$. 
For our applications to time series, the case $q_t^x=q\left(\cdot\vert x_t^{-}\right)$ will be of interest. The two following assumptions will be needed.

\begin{description}
\item
{\bf A1} The applications $(y,z,x)\mapsto q^{x}_t\left(y\vert z\right)$ and $(y,z,x)\mapsto q^{x}_t\left(y\vert z\right)$ are measurable and take positive values.
\item
{\bf A2} Setting
$$b_m:=\sup_{t\in\Z}\sup_{x\in \mathcal{X}}\sup_{y\stackrel{m}=\overline{y}}d_{TV}\left(q^{x}_t(\cdot\vert y),q^{x}_t(\cdot\vert \overline{y})\right),$$
we have $b_0<1$ and $\lim_{m\rightarrow \infty}b_m=0$.
\end{description}

Let us now introduce some additional notations. In what follows, we fix $t_0\in\Z$. For $z\in E^{\N}$ and $x\in \mathcal{X}$,
we denote by $Q_{t_0,x,z}$ 
the probability distribution on $\left(E^{\N},\mathcal{P}(E)^{\otimes \N}\right)$
defined by 
$$Q_{t_0,x,z}\left(\prod_{i=1}^n\{y_i\}\times \prod_{i=n+1}^{\infty}E\right)=\prod_{i=1}^nq_{t_0+i}^{x}\left(y_i\vert y_{i-1}^{-}\right)$$
with the convention $y_{-j}=z_j$ for $j\geq 0$.
We define $\overline{Q}_{t_0,x,z}$ in the same way, replacing the transition kernels $q_t^{x}$ with $\overline{q}_t^{x}$ in the previous expression.

\begin{lem}\label{central}
Assume that Assumptions {\bf A1-A2} hold true. Then for any $x\in\mathcal{X}$ and any couple $(z,\overline{z})\in E^{\N}\times E^{\N}$, there exists a probability measure $\widetilde{Q}_{t_0,x,z,\overline{z}}$ on $\left(E^{\N^{*}}\times E^{\N^{*}},\mathcal{P}(E)^{\otimes \N^{*}}\otimes\mathcal{P}(E)^{\otimes \N^{*}}\right)$ such that the three following conditions are satisfied. 
\begin{enumerate}
\item For $A,B\in \mathcal{P}(E)^{\otimes \N^{*}}$, we have
\begin{equation}\label{dyn0}
\widetilde{Q}_{t_0,x,z,\overline{z}}\left(A\times E^{\N^{*}}\right)=Q_{t_0,x,z}(A),\quad \widetilde{Q}_{t_0,x,z,\overline{z}}\left(E^{\N^{*}}\times B\right)=Q_{t_0,x,\overline{z}}(B).
\end{equation}
\item
For $t\geq 1$,
\begin{eqnarray}\label{dyn1}
&&\widetilde{Q}_{t_0,x,z,\overline{z}}\left(\left\{\left(y,\overline{y}\right)\in E^{\N^{*}}\times E^{\N^{*}}: y_t\neq \overline{y}_t\right\}\right)
\nonumber\\
&\leq& b_{t-1}^{*}+\sup_{s\in E^{\N}}d_{TV}\left(q^{x}_{t+t_0}(\cdot\vert s),\overline{q}^{x}_{t+t_0}(\cdot\vert s)\right)\nonumber\\
&+&\sum_{\ell=0}^{t-2}b_{\ell}^{*}\sup_{s\in E^{\N}}d_{TV}\left(q^{x}_{t+t_0-\ell-1}(\cdot\vert s),\overline{q}^{x}_{t+t_0-\ell-1}(\cdot\vert s)\right),
\end{eqnarray}
where $b_0^{*}=b_0$ and for $n\geq 1$, $b_n^{*}$ is equal to $\P\left(S_n^{(b)}=0\right)$ where $\left(S_n^{(b)}\right)_{n\geq 0}$ is a time-homogeneous Markov chain, starting at $0$ and with transition matrix $P$ defined by 
$$P(i,i+1)=1-b_i,\quad P(i,0)=b_i,\quad i\in\N.$$
\item
For all $C\in \mathcal{P}(E)^{\otimes \N^{*}}\otimes \mathcal{P}(E)^{\otimes \N^{*}}$, the application $(x,z,\overline{z})\mapsto \widetilde{Q}_{t_0,x,z,\overline{z}}(C)$ is measurable as an application from $\left(\mathcal{X}\times E^{\N}\times E^{\N},\mathcal{B}(\mathcal{X})\otimes \mathcal{P}(E)^{\otimes \N}\otimes \mathcal{P}(E)^{\otimes \N}\right)$ to $\left([0,1],\mathcal{B}([0,1])\right)$.
\end{enumerate}
\end{lem}

\paragraph{Notes} 
\begin{enumerate}
\item
Lemma \ref{central} is a central result for getting an upper bound of the 
total variation distance between the finite-dimensional distributions of two chains with complete connections and satisfying (\ref{dyn0}). For $i\in\N^{*}$, we denote by 
$y_i$ (resp. $\overline{y}_i$) the application from $E^{\N^{*}}\times E^{\N^{*}}$ to $E$ defined by $y_i(w,\overline{w})=w_i$ (resp. $\overline{y}_i(w,\overline{w})=\overline{w}_i$), $w,\overline{w}\in E^{\N^{*}}$.
Let $1\leq s\leq\ell$. If $Q_{t_0,x,z}^{(s,\ell)}$ and $\overline{Q}^{(s,\ell)}_{t_0,x,z}$ denote the
restriction of $Q_{t_0,x,z}$ (resp. $\overline{Q}_{t_0,x,z}$) to $\sigma(y_i: s\leq i\leq \ell)$, we have
\begin{eqnarray}\label{total}
d_{TV}\left(Q_{t_0,x,z}^{(s,\ell)},\overline{Q}_{t_0,x,\overline{z}}^{(s,\ell)}\right)&\leq& \widetilde{Q}_{t_0,x,z,\overline{z}}\left(y_t\neq \overline{y}_t; s\leq t\leq \ell\right)\nonumber\\ 
&\leq& \sum_{t=s}^{\ell}\widetilde{Q}_{t_0,x,z,\overline{z}}\left(y_t\neq \overline{y}_t\right)
\end{eqnarray}
and the total variation distance can be then bounded from (\ref{dyn1}).
\item
When the $q_t^{x}\equiv\overline{q}_t^{x}\equiv q$ and setting $\widetilde{Q}_{t_0,x,z,\overline{z}}=\widetilde{Q}_{t_0,z,\overline{z}}$, Lemma \ref{central} shows that 
$$\widetilde{Q}_{t_0,z,\overline{z}}\left(y_t\neq \overline{y}_t\right)\leq b^{*}_{t-t_0-1}$$
and we simply get control the total variation distance between the marginals at time $t$, when a time-homogeneous chain with complete connections is initialized with two different sequences. 
Such result has been proved by \citet{Bres} under a log-continuity assumption for the transition kernel $q$. Since we use an assumption slightly weaker in {\bf S2}, we will rewrite a detailed proof for the previous bound using our assumptions.
\item
One can note that the control of the total variation distance in Assumption {\bf A2} is uniform with respect to $x,t$. Getting a similar result without this uniformity is challenging but we did not find a way to relax it.  
\end{enumerate}

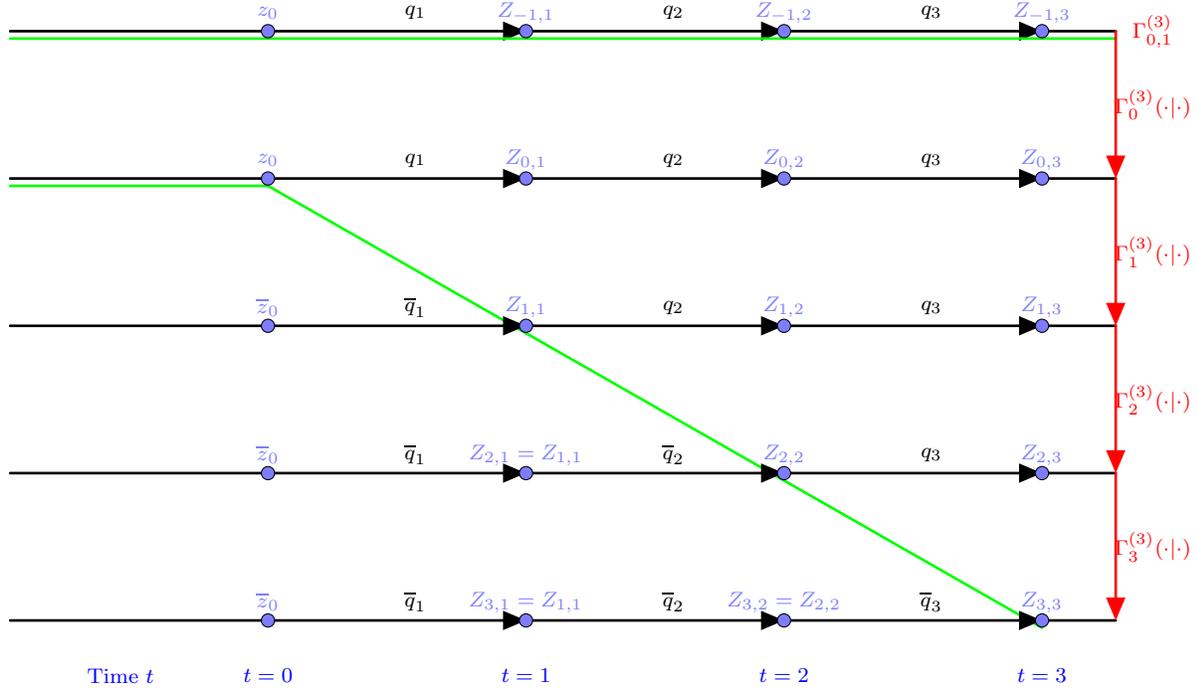
\begin{figure}
\begin{center}
\definecolor{ffqqqq}{rgb}{1,0,0}\definecolor{ududff}{rgb}{0.30196078431372547,0.30196078431372547,1}\definecolor{xdxdff}{rgb}{0.49019607843137253,0.49019607843137253,1}
\begin{tikzpicture}[line cap=round,line join=round,>=triangle 45,x=0.98cm,y=0.98cm]\clip(-15.5,-5) rectangle (7.661920883513466,5.311051711229949);
\draw [ line width=1pt] (-15.5,4) -- (-12,4);
\draw [->,line width=1pt] (-12,4) -- (-8.5,4);
\draw [->,line width=1pt] (-8.5,4) -- (-5,4);
\draw [->,line width=1pt] (-5,4) -- (-1.5,4);
\draw [line width=1pt] (-1.5,4) -- (-0.5,4);
\draw [green,line width=1pt] (-15.5,3.9) -- (-12,3.9);
\draw [green,line width=1pt] (-12,3.9) -- (-8.5,3.9);
\draw [green,line width=1pt] (-8.5,3.9) -- (-5,3.9);
\draw [green,line width=1pt] (-5,3.9) -- (-1.5,3.9);
\draw [green,line width=1pt] (-1.5,3.9) -- (-0.5,3.9);
\begin{scriptsize}
\draw [fill=xdxdff] (-12,4) circle (2.5pt);
\draw [fill=xdxdff] (-8.5,4) circle (2.5pt);
\draw [fill=xdxdff] (-5,4) circle (2.5pt);
\draw [fill=xdxdff] (-1.5,4) circle (2.5pt);
%\draw[color=black] (-14,4.25) node {$\mu$};
\draw[color=xdxdff] (-12,4.25) node {$z_0$};
\draw[color=xdxdff] (-8.5,4.25) node {$Z_{-1,1}$};
\draw[color=black] (-10,4.25) node {$q_1$};
\draw[color=xdxdff] (-5,4.25) node {$Z_{-1,2}$};
\draw[color=black] (-6.5,4.25) node {$q_2$};
\draw[color=xdxdff] (-1.5,4.25) node {$Z_{-1,3}$};
\draw[color=black] (-3,4.25) node {$q_3$};
\end{scriptsize}
\draw [green,line width=1pt] (-15.5,1.9) -- (-12,1.9);
\draw [green,line width=1pt] (-12,1.9) -- (-8.5,-0.1);
\draw [green,line width=1pt] (-8.5,-0.1) -- (-5,-2.1);
\draw [green,line width=1pt] (-5,-2.1) -- (-1.5,-4.1);

\draw [line width=1pt] (-15.5,2) -- (-12,2);
\draw [->,line width=1pt] (-12,2) -- (-8.5,2);
\draw [->,line width=1pt] (-8.5,2) -- (-5,2);
\draw [->,line width=1pt] (-5,2) -- (-1.5,2);
\draw [line width=1pt] (-1.5,2) -- (-0.5,2);
\begin{scriptsize}
\draw [fill=xdxdff] (-12,2) circle (2.5pt);
\draw [fill=xdxdff] (-8.5,2) circle (2.5pt);
\draw [fill=xdxdff] (-5,2) circle (2.5pt);
\draw [fill=xdxdff] (-1.5,2) circle (2.5pt);
%\draw[color=black] (-14,2.25) node {$\overline{\mu}$};
\draw[color=xdxdff] (-12,2.25) node {$z_0$};
\draw[color=xdxdff] (-8.5,2.25) node {$Z_{0,1}$};
\draw[color=black] (-10,2.25) node {$q_1$};
\draw[color=xdxdff] (-5,2.25) node {$Z_{0,2}$};
\draw[color=black] (-6.5,2.25) node {$q_2$};
\draw[color=xdxdff] (-1.5,2.25) node {$Z_{0,3}$};
\draw[color=black] (-3,2.25) node {$q_3$};
\end{scriptsize}

\draw [line width=1pt] (-15.5,0) -- (-12,0);
\draw [->,line width=1pt] (-12,0) -- (-8.5,0);
\draw [->,line width=1pt] (-8.5,0) -- (-5,0);
\draw [->,line width=1pt] (-5,0) -- (-1.5,0);
\draw [line width=1pt] (-1.5,0) -- (-0.5,0);
\begin{scriptsize}
\draw [fill=xdxdff] (-12,0) circle (2.5pt);
\draw [fill=xdxdff] (-8.5,0) circle (2.5pt);
\draw [fill=xdxdff] (-5,0) circle (2.5pt);
\draw [fill=xdxdff] (-1.5,0) circle (2.5pt);
%\draw[color=black] (-14,0.25) node {$\overline{\mu}$};
\draw[color=xdxdff] (-12,0.25) node {$\overline{z}_0$};
\draw[color=xdxdff] (-8.5,0.25) node {$Z_{1,1}$};
\draw[color=black] (-10,0.25) node {$\overline{q}_1$};
\draw[color=xdxdff] (-5,0.25) node {$Z_{1,2}$};
\draw[color=black] (-6.5,0.25) node {$q_2$};
\draw[color=xdxdff] (-1.5,0.25) node {$Z_{1,3}$};
\draw[color=black] (-3,0.25) node {$q_3$};
\end{scriptsize}

\draw [line width=1pt] (-15.5,-2) -- (-12,-2);
\draw [->,line width=1pt] (-12,-2) -- (-8.5,-2);
\draw [->,line width=1pt] (-8.5,-2) -- (-5,-2);
\draw [->,line width=1pt] (-5,-2) -- (-1.5,-2);
\draw [line width=1pt] (-1.5,-2) -- (-0.5,-2);
\begin{scriptsize}
\draw [fill=xdxdff] (-12,-2) circle (2.5pt);
\draw [fill=xdxdff] (-8.5,-2) circle (2.5pt);
\draw [fill=xdxdff] (-5,-2) circle (2.5pt);
\draw [fill=xdxdff] (-1.5,-2) circle (2.5pt);
%\draw[color=black] (-14,-1.75) node {$\overline{\mu}$};
\draw[color=xdxdff] (-12,-1.75) node {$\overline{z}_0$};
\draw[color=xdxdff] (-8.5,-1.75) node {$Z_{2,1}=Z_{1,1}$};
\draw[color=black] (-10,-1.75) node {$\overline{q}_1$};
\draw[color=xdxdff] (-5,-1.75) node {$Z_{2,2}$};
\draw[color=black] (-6.5,-1.75) node {$\overline{q}_2$};
\draw[color=xdxdff] (-1.5,-1.75) node {$Z_{2,3}$};
\draw[color=black] (-3,-1.75) node {$q_3$};
\end{scriptsize}

\draw [line width=1pt] (-15.5,-4) -- (-12,-4);
\draw [->,line width=1pt] (-12,-4) -- (-8.5,-4);
\draw [->,line width=1pt] (-8.5,-4) -- (-5,-4);
\draw [->,line width=1pt] (-5,-4) -- (-1.5,-4);
\draw [line width=1pt] (-1.5,-4) -- (-0.5,-4);
\begin{scriptsize}
\draw [fill=xdxdff] (-12,-4) circle (2.5pt);
\draw [fill=xdxdff] (-8.5,-4) circle (2.5pt);
\draw [fill=xdxdff] (-5,-4) circle (2.5pt);
\draw [fill=xdxdff] (-1.5,-4) circle (2.5pt);
%\draw[color=black] (-14,-3.75) node {$\overline{\mu}$};
\draw[color=xdxdff] (-12,-3.75) node {$\overline{z}_0$};
\draw[color=xdxdff] (-8.5,-3.75) node {$Z_{3,1}=Z_{1,1}$};
\draw[color=black] (-10,-3.75) node {$\overline{q}_1$};
\draw[color=xdxdff] (-5,-3.75) node {$Z_{3,2}=Z_{2,2}$};
\draw[color=black] (-6.5,-3.75) node {$\overline{q}_2$};
\draw[color=xdxdff] (-1.5,-3.75) node {$Z_{3,3}$};
\draw[color=black] (-3,-3.75) node {$\overline{q}_3$};
\draw[color=blue](-14,-4.75) node {Time $t$};
\draw[color=blue](-12,-4.75) node {$t=0$};
\draw[color=blue](-8.5,-4.75) node {$t=1$};
\draw[color=blue](-5,-4.75) node {$t=2$};
\draw[color=blue](-1.5,-4.75) node {$t=3$};

\end{scriptsize}

\draw[<-,line width=1, color=red](-0.5,-4) -- (-0.5,-2);
\draw[<-,line width=1, color=red](-0.5,-2) -- (-0.5,0);
\draw[<-,line width=1, color=red](-0.5,0) -- (-0.5,2);
\draw[<-,line width=1, color=red](-0.5,2) -- (-0.5,4);
\begin{scriptsize}
\draw[color=red] (0,4) node {$\Gamma^{(3)}_{0,1}$};
\draw[color=red] (0,3) node {$\Gamma_0^{(3)}(\cdot\vert\cdot)$};
\draw[color=red] (0,1) node {$\Gamma_1^{(3)}(\cdot\vert\cdot)$};
\draw[color=red] (0,-1) node {$\Gamma_2^{(3)}(\cdot\vert\cdot)$};
\draw[color=red] (0,-3) node {$\Gamma_3^{(3)}(\cdot\vert\cdot)$};
\end{scriptsize}

%\draw [line width=1.3pt,domain=-17.52505:7.661920883513466] plot(\x,{(-7.934778104306005-0.6630112299465223*\x)/0.005339163763064647});
%\draw[color=ffqqqq] (1.145779461559198,3.587600828877008) node {$Sequential$};\draw[color=ffqqqq] (1.5,2) node {$Sequential_1$};
\end{tikzpicture}
\end{center}
\caption{Description of the coupling when $k=3$. $\Gamma_i$, which denotes the coupling between the path $i-1$ and the path $i$, is obtained by applying  iteratively the maximal coupling over the time index $t$. For $i=0$, only the distribution of past values ($t\leq 0$) is changed ($(Z_{-1,j})_{j\leq 0}$ and $(Z_{0,j})_{j\leq 0}$ are assumed to be independent) whereas for $i=1,2,3$, only one conditional distribution is modified between two successive paths.
The gluing technique allows to define all the paths on the same probability space by using the conditional distribution of the coupling measures. The two green lines correspond to the two paths $\left(y_t\right)$ and $(\overline{y}_t)$. \label{coupling}} 
\end{figure}

\paragraph{Proof of Lemma \ref{central}}
Without loss of generality, we will assume that $t_0=0$, the general case will follow by replacing $t$ by $t-t_0$ in the bound we will derive. We then remove the index $t_0$ from all our notations.
We will apply the technique of maximal coupling already used by \citet{Bres} for getting a bound on the relaxation speed of chains with complete connections. We defer the reader to \cite{Bres}, equation $4.9$ for a precise definition of the maximal coupling of two probability measures $\alpha$ and $\overline{\alpha}$ on the finite set $E$. In what follows, we will simply use the fact that there exists a probability measure $\alpha\widetilde{\times}\overline{\alpha}$ on $E\times E$ such that 
$$\alpha\widetilde{\times}\overline{\alpha}\left(\{(y,\overline{y})\in E^2: y\neq \overline{y}\}\right)=d_{TV}\left(\alpha,\overline{\alpha}\right).$$
For a sequence $\omega=(\omega_{i,j})_{i,j\geq 0}\in E^{\N\times \N^{*}}$, we denote, for $(j,k)\in\N\times\N^{*}$, by $\omega_{j,1:k}$ the vector $(\omega_{j,1},\ldots,\omega_{j,k})\in E^k$.
We then set 
$$\Gamma^{(k)}_{0,1}\left(\omega_{0,1:k}\right)=\prod_{t=1}^kq_t^{x}\left(\omega_{0,t}\vert \omega_{0,t-1}^{-}\right).$$
In the previous expressions and the next ones, we always use the convention $\omega_{0,-i}=z_i$ and $\omega_{j,-i}=\overline{z}_i$ for $i\geq 0$ and $j\geq 1$. $\Gamma^{(k)}_{0,1}$ defines a probability measure on $E^k$.
Next, we define $k$ probability kernels $\Gamma^{(k)}_1,\ldots,\Gamma^{(k)}_k$ on $E^k$ in the following way. 
$$\Gamma^{(k)}_1\left(\omega_{1,1:k}\vert \omega_{0,1:k}\right)=\frac{\prod_{t=1}^k\left[q_t^{x}(\cdot\vert\omega_{0,t-1}^{-})\widetilde{\times}q^{x}_t(\cdot\vert \omega_{1,t-1}^-)\right](\omega_{0,t},\omega_{1,t})}{\Gamma^{(k)}_{0,1}\left(\omega_{0,1:k}\right)}.$$
If $2\leq j\leq k$, the kernel $\Gamma^{(k)}_j$ is defined by the equality
\begin{eqnarray*}
&&\Gamma^{(k)}_j\left(\omega_{j+1,1:k}\vert \omega_{j,1:k}\right)\times\prod_{t=1}^{j-1}\overline{q}_t^{x}\left(\omega_{j,t}\vert \omega_{j,t-1}^{-}\right)\cdot \prod_{t=j}^kq_t^{x}\left(\omega_{j,t}\vert \omega_{j,t-1}^{-}\right) \\
&=&\prod_{t=1}^{j-1}\left[\overline{q}_t^{x}(\cdot\vert\omega_{j,t-1}^{-})\widetilde{\times}\overline{q}^{x}_t(\cdot\vert \omega_{j+1,t-1}^-)\right](\omega_{j,t},\omega_{j+1,t})\times \left[q_j^{x}(\cdot\vert\omega_{j,j-1}^{-})\widetilde{\times}\overline{q}^{x}_j(\cdot\vert \omega_{j+1,j-1}^-)\right](\omega_{j,j},\omega_{j+1,j})\\
&\times& \prod_{t=j+1}^k \left[q_t^{x}(\cdot\vert\omega_{j,t-1}^{-})\widetilde{\times}q^{x}_t(\cdot\vert \omega_{j+1,t-1}^-)\right](\omega_{j,t},\omega_{j+1,t}).
\end{eqnarray*}
Finally, we define a probability measure $P^{(k)}_{x,z,\overline{z}}$ on $(E^k)^{k+1}$ by 
$$P^{(k)}_{x,z,\overline{z}}\left(\omega_{0,1:k},\ldots,\omega_{k+1,1:k}\right)=\Gamma_{0,1}^{(k)}(\omega_{0,1:k}\times\prod_{j=0}^k\Gamma_j^{(k)}\left(\omega_{j+1,1:k}\vert \omega_{j,1:k}\right).$$
Let us give an interpretation of the measure $P^{(k)}_{x,z,\overline{z}}$. This measure is the probability distribution of the $k+1$ first coordinates of a Markov chain on the state space $E^k$. Each coordinate of the chain can be seen as a path of a chain with complete connection.
\begin{itemize}
\item
$\Gamma_{0,1}^{(k)}$ is the distribution of $k$ successive coordinates of a chain with complete connection with initialization $z_{-i}$ for $i\leq 0$ and transition kernels $q_1^{x},\ldots,q_k^{x}$.
\item
The joint distribution of the first path and the second path is obtained by applying iteratively the maximal coupling to the transition kernels $(q_t^{x},q_t^{x})$ from time $t=1$ to time $t=k$. The second path is initialized with $\overline{z}_{-i}$ for $i\leq 0$. The second path  has then the same transition kernels as the first path but a different initialization.   
\item
For $1\leq j\leq k$, the $j-$th path is initialized with $\overline{z}_{-i}$, $i\leq 0$ and has transition kernels $\overline{q}^{x}_1,\ldots, \overline{q}^{x}_{j-1},q^{x}_j,\ldots,q^{x}_k$. The path $j+1$ is obtained as the path $j$, except that at time $t=j$, the kernel $q_t^{x}$ is replaced with the kernel $\overline{q}_t^{x}$. The joint probability distribution of the paths $j$ and $j+1$ is obtained by applying iteratively the maximal coupling to these transition kernels. 
\end{itemize}
Our approach is equivalent to make several couplings of two successive paths having either a different  initialization or one transition kernel changing across the time and then "gluing" all the paths to define a joint probability distribution on $(E^k)^{k+1}$. Our definition of this joint probability measure is classical in coupling theory and can be seen as a particular application of the so-called gluing lemma. See \citet{Villani}, Chapter $1$.
It is much easier to visualize such a coupling graphically. Figure \ref{coupling} gives a description of this coupling scheme when $k=3$.

Next, let us observe that $(x,z,\overline{z})\mapsto P^{(k)}_{x,z,\overline{z}}$ is measurable.
This is a consequence of the definition of $P^{(k)}_{x,z,\overline{z}}$ and of the explicit expression of the maximal coupling of two discrete probability measures in term of the marginals.
Measurability of the previous application then follows from Assumption {\bf A1}. 
We now mention that the sequence $\left(P^{(k)}_{x,z,\overline{z}}\right)_{k\geq 1}$ satisfies Kolmogorov's compatibility conditions. Indeed, one can show that
$$P^{(k)}_{x,z,\overline{z}}\left(\omega_{0,1:k},\ldots,\omega_{k+1,1:k}\right)=\sum_{\omega_{k+2,1:k+1}\in E^{k+1}}\sum_{\omega_{0:k+1,k+1}\in E^{k+2}}P^{(k+1)}_{x,z,\overline{z}}\left(\omega_{0,1:k+1},\ldots,\omega_{k+2,1:k+1}\right).$$
From the Kolmogorov extension theorem, there exists a unique probability measure $P_{x,z,\overline{z}}$ on $E^{\N\times \N^{*}}$ compatible with this sequence. Note that, for any $A\in \mathcal{P}(E)^{\N\times \N^{*}}$, the application $(x,z,\overline{z})\mapsto P_{x,z,\overline{z}}(A)$ is still measurable.
This was already justified when $A$ is a cylinder set. Extension of the measurability for $A$ arbitrary follows from a monotone class argument.

Now for $\omega\in E^{\N\times \N^{*}}$ and $j\geq -1$, $k\geq 1$, we define $Z_{j,k}(\omega)=\omega_{j+1,k}$. We define the probability distribution $\widetilde{Q}_{x,z,\overline{z}}$ as the pushforward measure of $P_{x,z,\overline{z}}$ obtained from $\left((Z_{-1,t})_{t\geq 1},(Z_{t,t})_{t\geq 1}\right)$.
Note that, from our construction with the maximal coupling, we have automatically $Z_{k,t}=Z_{t,t}$ for $k\geq t\geq 0$, $P_{x,z,\overline{z}}$ a.s.
Indeed, when the two past sequences are equal, the maximal coupling generates two identical random variables. We then deduce that $(Z_{t,t})_{t\geq 1}$ has transition kernels $(\overline{q}_t^{x})_{t\geq 1}$. This proves (\ref{dyn0}).
Let us now prove the bound (\ref{dyn1}). Let $t$ be a positive integer. 
We denote by $E_{x,z,\overline{z}}$ the mathematical expectation under $P_{x,z,\overline{z}}$.
From the triangular inequality, we have
\begin{equation}\label{montre}
\widetilde{Q}_{x,z,\overline{z}}\left(\{y_t\neq \overline{y}_t\}\right)=P_{x,z,\overline{z}}\left(Z_{-1,t}\neq Z_{t,t}\right)\leq \sum_{k=-1}^{t-1}P_{x,z,\overline{z}}\left(Z_{k,t}\neq Z_{k+1,t}\right).
\end{equation}
If $k=-1$, one can use Lemma \ref{auxilliaire} to get
$$P_{x,z,\overline{z}}\left(Z_{-1,t}\neq Z_{0,t}\right)\leq b_{t-1}^{*}.$$
Indeed $Z_{-1,\cdot}$ and $Z_{0,\cdot}$ are constructed using the maximal coupling and when $t\geq 1$, 
the transition kernel for the two paths equals $q_t(\cdot\vert\cdot)$.
If $t=k+1\geq 1$, we have from the definition of the maximal coupling and from our construction
\begin{eqnarray*}
P_{x,z,\overline{z}}\left(Z_{k,t}\neq Z_{k+1,t}\right)&=& E_{x,z,\overline{z}}\left[P_{x,z,\overline{z}}\left(Z_{k,t}\neq Z_{k+1,t}\vert \sigma\left(Z_{k,t-j},Z_{k+1,t-j}; j\geq 1\right)\right)\right]\\
&\leq& \sup_{w\in E^{\N}}d_{TV}\left(q^{x}_t\left(\cdot\vert w\right),\overline{q}^{x}_t\left(\cdot\vert w\right)\right).
\end{eqnarray*}
Next, if $t\geq k+2\geq 2$, we have
\begin{eqnarray*}
&&P_{x,z,\overline{z}}\left(Z_{k,t}\neq Z_{k+1,t}\right)\\&=&
E_{x,z,\overline{z}}\left[P_{x,z,\overline{z}}\left(Z_{k,t}\neq Z_{k+1,t}\vert\sigma\left(Z_{k,k+1-j},Z_{k+1,k+1-j}; j\geq 0\right)\right)\mathds{1}_{Z_{k,k+1}\neq Z_{k+1,k+1}}\right]\\
&\leq& b^{*}_{t-k-2}P_{x,z,\overline{z}}\left(Z_{k,k+1}\neq Z_{k+1,k+1}\right)\\
&\leq&  b^{*}_{t-k-2}\sup_{w\in E^{\N}}d_{TV}\left(q^{x}_{k+1}\left(\cdot\vert w\right),\overline{q}^{x}_{k+1}\left(\cdot\vert w\right)\right).
\end{eqnarray*}
Let us comment the previous bounds. The first equality follows from the fact that on the event $\left\{Z_{k,k+1}=Z_{k+1,k+1}\right\}$, we automatically have $Z_{k,j}=Z_{k+1,j}$ for $j\geq k+1$. This is due to the maximal coupling and to the fact that, from our construction, we have $Z_{k,j}=Z_{k+1,j}$ for $j\leq k$. The second bound follows from Lemma \ref{auxilliaire} and the fact that for $s=k+2,\ldots,t$,
the two paths $Z_{k,\cdot}$ and $Z_{k+1,\cdot}$ have the same transition kernels $q^{x}_s(\cdot\vert\cdot)$, $s\geq 1$. Finally, the third bound follows from the definition of the maximal coupling.
The bound (\ref{dyn1}) follows from (\ref{montre}) and our previous bounds.
Finally, the third point of the lemma follows from the measurability properties of $(x,z,\overline{z})\mapsto P_{x,z,\overline{z}}$. The proof of Lemma \ref{central} is now complete.$\square$

Next we provide a perturbation result for homogeneous chains with complete connections, i.e. $q_t\equiv q$ and $\overline{q}_t\equiv\overline{q}$ for all integer $t$. This result will not be used in the rest of the paper. However, it extends a standard perturbation result
for finite-state Markov chains and has then an independent interest.
If $\sum_{m\geq 1}b_m<\infty$ then $\sum_{m\geq 1}b_m^{*}<\infty$ (see \citet{Bres}, Proposition $2$). In this case, $b_m^{*}\rightarrow 0$ and there exists a unique stationary chain $(Y_k)_{k\in\Z}$ with complete connection and transition kernel $q$. Existence and unicity can hold under weaker conditions. See \citet{Bres}, Remark $1$.
We will also assume that the transition kernel $\overline{q}$ satisfies Assumptions ${\bf A1-A2}$ with summable coefficients $\overline{b}_m$.
By setting $t_0=0$ and letting $t$ going to $-\infty$ in Lemma \ref{central}, we obtain the following result.

\begin{cor}\label{pert}
Assume that $\sum_{m\geq 1}(b_m+\overline{b}_m)<\infty$ and let $\pi$ ($\overline{\pi}$ resp.) be the marginal distribution of the chain with transition kernel $q$ ($\overline{q}$ resp.). Then
$$d_{TV}\left(\pi,\overline{\pi}\right)\leq \left(1+\sum_{m\geq 0}b^{*}_m\right)\cdot \sup_{y\in E^{\N}}d_{TV}\left(q\left(\cdot\vert y\right),\overline{q}\left(\cdot\vert y\right)\right).$$
\end{cor}

Corollary \ref{pert} shows that the marginal distribution of the chain is a Lipschitz functional of its transition kernel. 
Let us detail this result in the Markov case, i.e. $q(\cdot\vert\cdot)$ is a stochastic matrix on $E$. In this case $b_m=0$ for $m\geq 1$ and it is easily seen that $b_m^{*}=b_0^m$.
We obtain  $1+\sum_{m\geq 0}b_m^{*}=(1-b_0)^{-1}$. We then recover a basic result for the perturbation of Markov chain using the ergodicity coefficient $b_0$ of the Markov chain with transition $q$. See for instance \citet{MIT}, Theorem $3.2$.

\section{Stationary categorical time series models with covariates}\label{3}

In this section, we consider a finite set $E$ with cardinal $N$. 
We will consider a stationary covariate process $X=(X_t)_{t\in \Z}$ taking values in $(\R^d,\vert\cdot\vert)$ where $\vert\cdot\vert$ is a norm on $\R^d$. 
For a sequence $(x_t)_{t\in\Z}$ and $t\in\Z$, we use the notation $x_t^{-}=\left(x_{t-j}\right)_{j\geq 0}$.
Let $(Y_t)_{t\in\Z}$ a time series taking values in $E$ and such that
\begin{equation}\label{spec} 
\P\left(Y_t=w\vert Y_{t-1}^{-},X\right)=q\left(w\vert Y_{t-1}^{-},X_t^{-}\right),\quad t\in \Z.
\end{equation}
We assume that the applications $(w,y,x)\mapsto q\left(w\vert  y, x_t^{-}\right)$ are measurable, as  applications 
from $E\times E^{\N}\times \mathcal{D}$ to $(0,1)$, where $\mathcal{D}\in \mathcal{B}\left(\R^d\right)^{\otimes \Z}$ is such that $\P\left(X\in \mathcal{D}\right)=1$. Moreover, we impose $\sum_{w\in E}q\left(w\vert y, x\right)=1$ for all $(y,x)\in E^{\N}\times \mathcal{D}$.

\subsection{Existence of a stationary and ergodic solution}
The following assumptions will be needed.

\begin{description}
\item[S1] The covariate process $X=(X_t)_{t\in\Z}$ stationary and ergodic.

\item [S2]
Setting for $m\geq 0$,
$$b_m=\sup\left\{d_{TV}\left(q\left(\cdot\vert y,x_t^{-}\right),q\left(\cdot \vert y',x_t^{-}\right)\right): (y,y',x)\in E^{\N}\times E^{\N}\times \mathcal{D}, t\in\Z, y\stackrel{m}=y'\right\},$$ 
we have $b_0<1$ and $\sum_{m\geq 0}b_m<\infty$.
\end{description}

\paragraph{Note.} 
Assumption {\bf S2} guarantees that $\sum_{m\geq 0}b^{*}_m<\infty$, where the $b_m^{*}$'s are related to the $b_m'$s as described in Lemma \ref{central}. For a proof, see \citet{Bres}, Proposition $2$.
Basically, the decrease of the sequence $\left(b_m^{*}\right)_{m\geq 0}$ is of the same order as the sequence $(b_m)_{m\geq 0}$. 
One can note that we impose a control of the total variation distances which is uniform with respect to the path of the covariate process $X$. We did not find a solution for removing this assumption.
However, as we will see in the examples, when the contribution of the covariates is additive in some generalized linear models, this assumption is often satisfied even if the covariate process is unbounded.

\begin{theo}\label{mainstat}
Assume that the assumptions {\bf S1-S3} hold true. 
\begin{enumerate}
\item
There then exists a unique stochastic processes $(Y_t)_{t\in \Z}$ satisfying (\ref{spec}). Moreover for any bounded measurable function $h:E^{\N}\rightarrow \R$, we have 
$$\E\left[h\left(Y_t^{-}\right)\vert X\right]=\E\left[h\left(Y_t^{-}\right)\vert X_t^{-}\right].$$

\item
The bivariate process $\left((Y_t,X_t)\right)_{t\in\Z}$ is stationary and ergodic. 
\end{enumerate}
\end{theo}

\paragraph{Proof of Theorem \ref{mainstat}}
\begin{enumerate}
\item
Let $\left(\overline{\Omega},\overline{\mathcal{A}},\overline{\P}\right)$ be a probability space on which the covariate process $X$ is defined. For simplicity, we assume that $\overline{\Omega}=G^{\Z}$ is the canonical space of the paths. We then have $X_t(\omega)=\omega_t$ for all $(t,\omega)\in \Z\times\overline{\Omega}$. 
Existence of a stochastic process $(Y_t)_{t\in\Z}$ satisfying (\ref{spec}) is understood as follows. We consider an enlargement $\left(\Omega,\mathcal{A},\P\right)$ 
of the initial probability space with $\Omega=E^{\Z}\times \overline{\Omega}$, 
$\mathcal{A}=\mathcal{P}(E)^{\otimes\Z}\otimes \overline{\mathcal{A}}$. For all $(y,\overline{\omega},t)\in\Omega\times \Z$, we set $Y_t(y,\overline{\omega})=y$, $X'_t(y,\overline{\omega})=X_t(\overline{\omega})$ and $\P$ is
the probability measure such that $\P\left(X'\in E^{\Z}\times B\right)=\overline{\P}\left(X\in B\right)$ for all $B\in\overline{\mathcal{A}}$ and
for which (\ref{spec}) is satisfied, replacing $X$ with $X'$.
Unicity is understood as follows: if $\left(\Omega',\mathcal{A}',\P'\right)$ is another probability space on which two processes ${\bf X'}$ and ${\bf Y'}$ are defined and satisfy
$\P'\left(X'\in B\right)=\P\left(X\in B\right)$ and (\ref{spec}), then $\P'\left((Y',X')\in A\right)=\P\left((Y,X)\in A\right)$ for all $A\in \mathcal{P}(E)^{\otimes \Z}\otimes \mathcal{B}(G)^{\otimes\Z}$.

To show these properties, we will construct, for each $x\in\mathcal{D}$, a family of finite-dimensional probability distributions $\left\{\nu^{I}_{x}: I=\{s+1,\ldots,s+n\}; (s,n)\in \Z\times\N^{*}\right\}$
such that for all $I$, $\nu_{x}^{I}$ is a probability measure on $E^I$ and satisfies the compatibility conditions of the Kolmogorov's extension theorem.
To this end, we will use the coupling result of Lemma \ref{central} and more precisely the control in total variation given in (\ref{total}). 
Let us consider a set $I=\{s+1,\ldots,s+n\}$ of $n$ successive integers and two elements $z,\overline{z}$ of $E^{\N}$.   
For an integer $i\geq 2$, set $t_0=s-i$ and for $t\geq t_0+1$ and $y\in E^{\Z}$,
$$q^{x}_t(y_0\vert y_{-1}^{-})=\overline{q}^{x}_t(y_0\vert y_{-1}^{-})=q\left(y_0\vert y_{-1}^{-},x^{-}_t(\omega)\right).$$
We define a probability measure $\nu^{I}_{x,i,z}$ by 
$$\nu^{I}_{x,i,z}(y_{s+1},\ldots,y_{s+n})=\int\sum_{y_{t_0+1},\ldots,y_s}\prod_{t=t_0+1}^{s+n}q^{x}_t(y_t\vert y_{t-1}^{-})\delta_z\left(d y_{t_0}^{-}\right).$$
Using (\ref{total}), we have 
$$d_{TV}\left(\nu^{I}_{x,i,z},\nu^{I}_{x,i,\overline{z}}\right)\leq  
\sum_{\ell=0}^{n-1} b^{*}_{i+\ell}.$$
Assumption {\bf S2} guarantees the summability of the $b_m^{*}$'s and hence that $b^{*}_m\rightarrow 0$. Moreover, the previous bound does not depend on the couple $(z,\overline{z})$ and goes to zero when $i\rightarrow \infty$, one can show that the sequence $\left(\nu^{I}_{x,i,z}\right)_{i\geq 2}$ is a Cauchy sequence in the simplex of $\R^n$ and has a limit which does not depend on $z$. We then set 
\begin{equation}\label{etap}
\nu^{I}_{x}=\lim_{i\rightarrow \infty}\nu^{I}_{x,i,z}.
\end{equation}
The compatibility conditions on the family of finite dimensional distributions
$$\mathcal{G}=\left\{\nu^{I}_{x}: I=\{s+1,\ldots,s+n\}, (s,n)\in \Z\times\N^{*}\right\}$$
follows from the fact that for any $i\geq 2$ and $I=\{s+1,\ldots,s+n\}$,
$$\sum_{y_{s+n+1}\in E}\nu^{I\cup\{s+n+1\}}_{x,z,i}(y_{s+1},\ldots,y_{s+n+1})=\nu^{I}_{x,z,i}\left(y_{s+1},\ldots,y_{s+n}\right),$$
$$\sum_{y_s\in E}\nu^{I\cup\{s\}}_{x,z,i}(y_s,\ldots,y_{s+n})=\nu^{I}_{x,z,i+1}\left(y_{s+1},\ldots,y_{s+n}\right).$$
The Kolmogorov's extension theorem guarantees existence and unicity of a probability measure $\nu_{x}$ on $\left(\E^{\Z},\mathcal{P}(E)^{\Z}\right)$ compatible with the family $\mathcal{G}$.
We then $\Omega=E^{\Z}\times \overline{\Omega}$ and for $\omega=(y,\overline{\omega})\in \Omega$, $\P(d\omega)=\nu_{X(\overline{\omega})}(dy)\overline{\P}(d\overline{\omega})$ and $Y_t\left(\omega\right)=y_t$ for $t\in\Z$.
We point out that measurability of the application $x\mapsto \nu_x(A)$ can be shown first when $A$ is a cylinder set and then for an arbitrary $A\in \mathcal{P}(E)^{\otimes\Z}$ using a monotone class argument. 

\item
Next we show that the process $(Y_t)_{t\in\Z}$ defined in the previous point satisfies (\ref{spec}).
This is equivalent to show that the probability measure $\nu_{x}$ defined in the previous point is compatible with the sequence $(q^{x}_t)_t$. We keep the notations of the previous point. Let $\epsilon\in(0,1)$ and $
h:E\rightarrow \R$ and $g:E^k\rightarrow \R$ be some functions bounded by one. Let $\overline{y}$
be an arbitrary element of $E$.
From our assumptions, there exists and integer $k\geq 1$ such that $b_m\leq \epsilon$ if $m\geq k$. Set 
$$q^{x}_{m,t}(y_t\vert y_{t-1:t-m})=q^{x}_t\left(y_t\vert y_{t-1},\ldots,y_{t-
m},\overline{y},\overline{y},\ldots\right).$$
We also set $I=\{t-k,\ldots,t\}$, $I_m=\{t-m,\ldots,t-1\}$ and we choose $i>m$ large 
enough such that 
$$d_{TV}\left(\nu_{x}^I,\nu^{I}_{x,z,i}\right)+d_{TV}\left(\nu_{x}^{I_m},\nu^{I_m}_{x,z,i}
\right)\leq \epsilon.$$
We have
$$\left\vert \int h(y_t)g(y_{t
-1},\ldots,y_{t-k})\left[d\nu^{I}_{x}(y_{t-k},\ldots,y_t)-d\nu^{I}_{x,z,i}(y_{t-k},\ldots,y_t)\right]\right\vert\leq 
\epsilon.$$
Next we set $q^{x}_{m,t}h(y_{t-1},\ldots,y_{t-m})=\sum_{y_t\in E}h(y_t)q^{x}_
{m,t}(y_t\vert y_{t-1:t-m})$ and 
$$A_m=\int q^{x}_{m,t} h(y_{t-1},\ldots,y_{t-m})g(y_{t-1},\ldots,y_{t-k})d\nu^{I_m}_{x,z,i}\left(y_{t-m},\ldots,y_{t-1}\right).$$
We have 
$$\left\vert A_m-\int h(y_t)g(y_{t
-1},\ldots,y_{t-k})d\nu^{I}_{x,z,i}(y_{t-k},\ldots,y_t)\right\vert\leq 
b_m\leq \epsilon.$$
Moreover 
$$\left\vert A_m-\int q^{x}_{m,t}h(y_{t-1},\ldots,y_{t-m})g(y_{t-1},\ldots,y_{t-k})d\nu^{I_m}_{x}\right\vert\leq \epsilon.$$
Using the fact that, 
$$\left\vert q^{x}_{m,t} h(y_{t-1},\ldots,y_{t-m})-q^{x}_t h\left(y_{t-1}^{-}\right)\right\vert\leq b_m\leq \epsilon,$$
we get 
$$\left\vert \int h(y_t)g(y_{t-1},\ldots, y_{t-k})d\nu^{I}_{x}(dy_t,\ldots,dy_{t-k})-\int q^{x}_t h(y_{t-1}^{-})g(y_{t-1},\ldots,y_{t-k})d\nu_{x}(y)\right\vert\leq 4\epsilon.$$
This proves that 
$$\int h(y_t)g(y_{t-1},\ldots, y_{t-k})d\nu^{I}_{x}(y_t,\ldots,y_{t-k})=\int q^{x}_t h(y_{t-1}^{-})g(y_{t-1},\ldots,y_{t-k})d\nu_{x}(y).$$
From a monotone class argument, we obtain (\ref{spec}).

\item
The equality between the two conditional expectations in point $1$ of Theorem \ref{mainstat} is a consequence of the expression of $\nu_x$ and $q_t^{x}$ (which only depends of $x_t^{-}$).

\item
Let us now show that the process $(V_t)_{t\in\Z}$ defined by $V_t=(Y_t,X_t)$  is stationary. It should be noticed first that if $I_t=\left\{t+1,\ldots,t+n\right\}$ for $t\in\Z$, then from the definition of the finite-dimensional distributions, we have
$$\nu^{I_t}_{x}=\nu^{I_0}_{\tau^t x} \mbox{ a.s.}$$
where $\tau^t x=(x_{t+j})_{j\in\Z}$. 
We then get for a measurable and bounded function $h:E^n\times \mathcal{D}\rightarrow \R$,
\begin{eqnarray*}
\E h\left(Y_{t+1},\ldots,Y_{t+n}, \tau^t X\right)&=& \sum_{y_1,\ldots,y_n\in E}\int h(y_1,\ldots,y_n,\tau^tX(\overline{\omega}))\nu^{I_t}_{X(\overline{\omega})}(y_1,\ldots,y_n)d\overline{\P}(\overline{\omega})\\
&=& \sum_{y_1,\ldots,y_n\in E}\int h(y_1,\ldots,y_n,\tau^tX(\overline{\omega})\nu^{I_0}_{\tau^tX(\overline{\omega})}(y_1,\ldots,y_n)d\overline{\P}(\overline{\omega})\\
&=& \sum_{y_1,\ldots,y_n}\int h(y_1,\ldots,y_n,\omega)\nu^{I_0}_{X(\overline{\omega})}(y_1,\ldots,y_n)d\overline{\P}(\overline{\omega})\\
&=& \E h\left(Y_1,\ldots,Y_n,X\right).
\end{eqnarray*}
This shows the stationarity of the process $\left(V_t\right)_{t\in\Z}$.

\item
Next, let us show uniqueness. Let $\left((Y_t',X_t')\right)_{t\in\Z}$ be a stationary process
satisfying the same assumptions. Setting $z_i={Y'}_{t-i}^{-}$, we know that from (\ref{etap}), we have
$\P'$ a.s.,
\begin{eqnarray*}
\lim_{i\rightarrow \infty}\P'\left(Y'_{t+1}=y_1,\ldots,Y'_{t+n}=y_n\vert \sigma\left(X',Y'_{t-j}:j\geq i\right)\right)&=&\lim_{i\rightarrow \infty}\nu^{I_t}_{{X',i+1,z_i}}(y_1,\ldots,y_n)\\
&=& \mu^{I_t}_{X'}(y_1,\ldots,y_n).
\end{eqnarray*}
Hence for any measurable and bounded function $h$, we have 
\begin{eqnarray*}
\E' h\left(Y'_{t+1},\ldots,Y'_{t+n},X'\right)&=& \lim_{i\rightarrow \infty}\E'\left[\E'\left(h\left(Y'_{t+1},\ldots,Y'_{t+n},X'\right)\vert\sigma(X',Y'_{t-j}:j\geq i)\right)\right]\\
&=& \sum_{y_1,\ldots,y_n}\E'\left[\mu^{I_t}_{X'}(y_1,\ldots,y_n)h(y_1,\ldots,y_n,X')\right]\\
&=& \sum_{y_1,\ldots,y_n}\E\left[\mu^{I_t}_{X}(y_1,\ldots,y_n)h(y_1,\ldots,y_n,{\bf X})\right]\\
&=& \E h\left(Y_{t+1},\ldots,Y_{t+n},X\right).
\end{eqnarray*}
The second equality follows from the Lebesgue theorem.

\item
For $t\in\Z$, we remind that $V_t=(Y_t,X_t)$. We now prove the ergodicity property for the process $\left(V_t\right)_{t\in\Z}$.
To this end, we adapt the direct proof of \citet{Kif} who proved ergodic properties of some Markov chains in random environments. 
Set $\mu=\P_X$, the probability distribution of $X$ under $\P$.
Remind that the measure $\nu_x$ constructed in point $1.$ is the probability distribution of $Y$ given $X=x$ and we will denote by $E_x$ the corresponding mathematical expectation. 
We will first consider the measure $\nu_x^{(0)}$, the probability distribution of $(Y_t)_{t\geq 0}$ given that $X=x$ and show that the measure $d\gamma(y,x)=d\nu^{(0)}_x(y)d\mu(x)$ is ergodic for 
the operator $(\theta,\tau).(y,x)=(\theta y,\tau x)$ where for $y\in E^{\N}$, $\theta(y)=(y_{t+1})_{t\in\N}$ and $\tau$ has been already defined as the shift operator on $(\R^d)^{\Z}$. 
For $t\in \Z$, ($t\in \N$ resp.), we denote by $y_t$ the coordinate application from $E^{\Z}$ ($E^{\N}$resp.) to $E$, i.e. $y_t(z)=z_t$ for $z\in E^{\Z}$ ($z\in E^{\N}$ resp.). Let $B\in \mathcal{P}(E)^{\otimes \N}$, $k\in\N$, $n$ an integer greater than $k$ and $w_0,\ldots,w_k\in E$. We have
\begin{eqnarray*}
\nu^{(0)}_x\left(y_0=w_0,\ldots,y_k=w_k,y\in \theta^{-n} B\right)
&=&\nu_x\left(y_0=w_0,\ldots,y_k=w_k,(y_t)_{t\geq 0}\in \theta^{-n} B\right)\\
&=& E_x\left[\prod_{i=0}^k\mathds{1}_{y_i=w_i}\times\nu_x\left((y_t)_{t\geq 0}\in\theta^{-n} B\vert y_k^{-}\right)\right].
\end{eqnarray*}
Using Lemma \ref{central} and the control of the total variation distance mentioned in the point $1.$ of the Notes, we also have
\begin{eqnarray*}
\left\vert \nu_x\left((y_t)_{t\geq 0}\in\theta^{-n} B\vert y_k^{-}\right)-
\nu^{(0)}_x\left(\theta^{-n}B\right)\right\vert
&\leq& 2 \sup_{z,\overline{z}}d_{TV}\left(Q_{k,x,z}^{(n-k+1,\infty)},Q_{k,x,\overline{z}}^{(n-k+1,\infty)}\right)\\
&\leq& 2\sum_{i=1}^{\infty}b^{*}_{n-k+i-1}\stackrel{n\rightarrow \infty}{\rightarrow}0.
\end{eqnarray*}
Note also that
$\nu^{(0)}_x\left(\theta^{-n}B\right)=\nu^{(0)}_{\tau^n x}\left(B\right)$.
We then get 
\begin{equation}\label{mixx}
\lim_{n\rightarrow \infty}\sup_{B\in \mathcal{P}(E)^{\otimes \N}}\left\vert \nu^{(0)}_x\left(A\cap \tau^{-n}B\right)-\nu^{(0)}_x(A)\nu^{(0)}_{\tau^n}(B)\right\vert=0,
\end{equation}
when $A$ is a cylinder set. Using approximation by finite unions of disjoint cylinder sets, one can extend (\ref{mixx}) to an arbitrary Borel set $A\in\mathcal{P}(E)^{\otimes \N}$.
Now let $\mathcal{I}$ be an invariant set in $E^{\N}\times\mathcal{D}$, i.e. $(\theta,\tau)^{-1}\mathcal{I}=\mathcal{I}$. It remains to show that $\gamma(\mathcal{I})\in \{0,1\}$.
We already mentioned in point $4.$, the equality $\nu^{(0)}_x(\theta^{-1}A)=\nu^{(0)}_{\tau x}(A)$ when $A$ is a cylinder set. This equality can be extended to any Borel set $A$. 
If $\mathcal{I}^x=\left\{y\in E^{\N}: (y,x)\in\mathcal{I}\right\}$, we have $\theta^{-1}\mathcal{I}^x=\mathcal{I}^{\tau^{-1}x}$.
From (\ref{mixx}), we then deduce that 
\begin{eqnarray*}
\nu_x^{(0)}\left(\mathcal{I}^x\right)-\nu_x^{(0)}\left(\mathcal{I}^x\right)^2
&=& \nu_x^{(0)}\left(\mathcal{I}^x\cap \theta^{-n}\mathcal{I}^{\tau^n x}\right)-\nu_x^{(0)}\left(\mathcal{I}^x\right)\nu_{\tau^nx}^{(0)}\left(\mathcal{I}^{\tau^n x}\right)\\
&\rightarrow& 0.
\end{eqnarray*}
This shows that $f(x):=\nu_x^{(0)}\left(\mathcal{I}^x\right)\in\{0,1\}$ for all $x$. 
Since $f(\tau x)=f(x)$, ergodicity of $X$ entails that we $\mu(\{f=1\})\in \{0,1\}$ and then $\gamma\left(\mathcal{I}\right)\in\{0,1\}$.
This shows that $\gamma$ is ergodic for $(\theta,\tau)$ and in particular that the process $\left(V_t\right)_{t\in\N}$ is ergodic. But this also entails ergodic properties for the two-sided sequence $(V_t)_{t\in\Z}$, see for instance Theorem $31$ in \citet{DDM} for a proof.$\square$
\end{enumerate}

\section{Examples}\label{4}
We now provide many examples of categorical time series models satisfying our assumptions.
We study in particular some observation-driven models proposed in the literature, which are parsimonious and then interesting for applications in statistics.

\subsection{Generalized linear model for binary time series}

Here we assume that $E=\{0,1\}$. We consider the following binary time series model defined by
\begin{equation}\label{autoregres0}
\P\left(Y_t=1\vert Y_{t-1}^{-},X\right)=F\left(\mu_t\right),\quad \mu_t=\sum_{j=1}^{\infty}a_j Y_{t-j}+\gamma'X_t,
\end{equation}
where $F$ is a cumulative distribution function, $(a_j)_{j\geq 1}$ is a summable sequence of real numbers and $\gamma\in\R^d$.
Model (\ref{autoregres0}) extends the model considered by \citet{comets} which does not contain exogenous regressors.

\begin{prop}\label{Bin0}
Assume that $F$ is Lipschitz, positive everywhere, $\sum_{j\geq 1} j\vert a_j\vert<\infty$ and $X$ satisfies Assumption {\bf S1}. There then exists a unique stationary process $(Y_t)_{t\in\Z}$ satisfying (\ref{autoregres0}). Moreover the bivariate process $\left((Y_t,X_t)\right)_{t\in\Z}$ is stationary and ergodic.
\end{prop}

\paragraph{Proof of Proposition \ref{Bin0}}
The result is a consequence of Theorem \ref{mainstat}. First we have $b_m\leq L\sum_{j\geq m}\vert a_j\vert$ with $L$ the Lipschitz constant of $F$. Our assumptions entails summability of the $b_m'$s.
The crucial point is to check the condition $b_0<1$.
Since the first term in the argument of $F$ is bounded, condition $b_0<1$ will follow if we show that for any $c>0$, 
\begin{equation}\label{chaud}
\sup_{\vert y\vert\leq c, z\in \R}\left\vert F(y+z)-F(z)\right\vert<1.
\end{equation}
Note first that $F$ has a limit at $\pm\infty$. Hence $\displaystyle\sup_{\vert y\vert\leq c,\vert z\vert >M}F(y+z)<1/4$ if $M$ is large enough. For such $M$, we also have $0\leq \displaystyle\inf_{\vert y\vert\leq c,\vert z\vert\leq M}F(y+z)\leq \displaystyle\sup_{\vert y\vert\leq c,\vert z\vert\leq M}F(y+z)<1$.
We then automatically have (\ref{chaud}) and then $b_0<1$. 
$\square$

 Model (\ref{autoregres}) is of theoretical interest but in practice observation-driven models lead to parsimonious representations of such dynamic. Let us consider the following version.

\begin{equation}\label{autoregres} 
\P\left(Y_t=1\vert Y_{t-1}^{-},X\right)=F\left(\mu_t\right),\quad \mu_t=\sum_{j=1}^q \beta_j\mu_{t-j}+\sum_{k=1}^p\alpha_k Y_{t-k}+\gamma'X_t,
\end{equation}
where $F$ is a cumulative distribution function $\alpha_1,\ldots,\alpha_p,\beta_1,\ldots,\beta_q\in \R$ and $\gamma\in\R^d$. 
We get the following result.

\begin{prop}\label{Bin1}
Assume that $F$ is positive everywhere and Lipschitz and that the covariate process $X$ satisfies {\bf S1} and $\E\log_+\vert X_1\vert$ for some $s\in(0,1)$. Assume further that the roots of the polynomial
$$\mathcal{P}(z)=1-\sum_{j=1}^q\beta_jz^j$$
are outside the unit disc. There then exists a unique stationary solution to (\ref{autoregres}). Moreover, 
the process $\left((Y_t,X_t)\right)_{t\in\Z}$ is stationary and ergodic.
\end{prop}

\paragraph{Notes}
\begin{enumerate}
\item
A classical choice for $F$ is the Gaussian c.d.f. (probit model) or the logistic c.d.f. (logistic model). Model of type (\ref{autoregres}) have been proposed but without a theoretical support 
by \citet{kauppi}, \citet{rydberg} or \citet{russell} for analyzing price changes or predicting recessions. When $\beta_1=\cdots=\beta_q=0$, a theory for the dynamic probit model can be found in \citet{deJong2011} or in \citet{FT1} who studied more general Markov models specified conditionally to some covariates. When there is no covariates, stationarity conditions for model (\ref{autoregres}) are given in \citet{FT1}. 
Our results then extend these previous contributions and also give a theoretical basis to some models used in econometrics.   
\item
 It is also possible to consider models with interactions  
between past values of the response and covariates. However, in general, application of Theorem \ref{mainstat} requires boundedness of the process $(\mu_t)_{t\in\Z}$ in (\ref{autoregres}). For simplicity, let us assume  
that $d=1$ and that the process is (conditionally to $X$) a first-order time-inhomogeneous Markov chain (called a Markov chain with covariates in what follows) with
$\mu_t=g\left(Y_{t-1},X_t\right)$
, $g:E\times\R\rightarrow \R$ being a measurable function.  
Then 
$$b_0=\sup_{z\in\R}d_{TV}\left(q(\cdot\vert 1,z),q(\cdot\vert 0,z)\right)=\sup_{z\in\R}\left\vert F\left(g(1,z)\right)-F\left(g(0,z)\right)\right\vert.$$
Assumption {\bf S2} is valid provided the second link function $g$ is bounded. When $g$ is not bounded, 
Assumption {\bf S2} is still valid when for any $z\in\R$, $g(0,z)$ and $g(1,z)$ have the same sign but this restriction seems to be quite artificial.
More generally, if $\mu_t=g\left(Y_{t-1},\ldots,Y_{t-p},X_t\right)$ with $g:E^p\times \R^d\rightarrow \R$ is measurable and bounded, Assumption {\bf S2} is satisfied with $b_0<1$ and $b_m=0$ if $m\geq p$.
We point out that these results are less sharp than that of \citet{FT1}, where existence and uniqueness of a stationary and ergodic solution for a Markov chain with covariates was obtained without this boundedness assumption.
However, Theorem \ref{mainstat} is compatible with non Markov processes	and then observation-driven models which are more difficult to study.
\end{enumerate}

\paragraph{Proof of Proposition \ref{Bin1}}
Setting $\lambda_t=\left(\mu_t,\ldots,\mu_{t-q+1}\right)'$, any solution of the problem (\ref{autoregres})
satisfies the recursions $\lambda_t=A\lambda_{t-1}+b_t$ with
$$A=\begin{pmatrix} \beta_1&\ldots&\beta_q\\& I_{q-1}& 0_{q-1,1}\end{pmatrix},\quad b_t=\begin{pmatrix} \sum_{k=1}^p\alpha_k Y_{t-k}+\gamma'X_t\\O_{q-1,1}\end{pmatrix},$$
where $0_{q-1,1}$ is a column vector of $0$ and $I_{q-1}$ is the identity matrix of size $(q-1)\times(q-1)$. Our assumptions guaranty that the spectral radius of $A$ is less than $1$. For a given operator norm $\Vert\cdot\Vert$, there then exists $r\in\N^{*}$ such that $\kappa:=\Vert A^r\Vert<1$. 
 One can then apply Lemma \ref{observation-driven} to show that any stationary solution $(Y_t)_{t\in\Z}$ satisfying (\ref{autoregres}) is a chain with complete connections and 
such that $b_m=O\left(\kappa^{m/r}\right)$. To end the proof, one can apply Theorem \ref{mainstat}. One only need to check that $b_0<1$. We observe that $q(1\vert y_{t-1}^{-},x_t^{-})$
is of the form $F\left(\sum_{j=1}^{\infty}\eta_j y_{t-j}+\sum_{k=0}^{\infty}\delta_k'X_{t-k}\right)$ for some summable sequences $(\eta_j)_{j\geq 1}$ and $(\delta_j)_{j\geq 0}$.
Since the first term in the argument of $F$ is bounded, condition $b_0<1$ follows exactly 
as in the proof of Proposition \ref{Bin0}, using (\ref{chaud}). Theorem \ref{mainstat} entails the result.$\square$

We still consider the binary case as in \ref{autoregres} with a non linear $\lambda_t$ and with one lag for simplicity. For a function $g:\R\mapsto \R$, we assume that 
\begin{equation}\label{nonlin}
\lambda_t=g\left(\lambda_{t-1}\right)+\alpha Y_{t-1}+\gamma' X_t.
\end{equation}
Such type of model has been proposed by \citet{russell} for analyzing financial transactions prices setting $g(s)=\beta s-\alpha F(s)$. Note that with the last specification, if $\vert \beta\vert<1$, $\lambda_t$ writes as a linear combination of the martingale differences $Y_{t-j}-F\left(\lambda_{t-j}\right)$, $j\geq 1$. 

\begin{prop}\label{nonlinear}
Assume that $F$ is positive everywhere, Lipschitz and $g$ is Lipschitz with 
$$\vert g(s)-g(s')\vert\leq \kappa \vert s-s'\vert,\quad (s,s')\in\R^2$$
for some $\kappa\in(0,1)$. Assume further that the process $X$ satisfies Assumption {\bf S1} and $\E\log_+\vert X_0\vert<\infty$.
There exists a unique stationary process $(Y_t)_{t\in\Z}$ satisfying (\ref{nonlin}). Moreover, the process $\left((X_t,Y_t)\right)_{t\in\Z}$ is stationary and ergodic.
\end{prop}

\paragraph{Proof of Proposition \ref{nonlinear}}
As in the proof of Proposition \ref{Bin1}, we use Lemma \ref{observation-driven}
which shows that any solution of (\ref{nonlin}) is a chain with complete connections for which the coefficients $b_m$ decay geometrically fast. To show {\bf S2}, it remains to show the condition $b_0<1$. Lemma \ref{observation-driven} shows that $\sup_{t,x,y,\overline{y}}\left\vert \lambda_t^{y,x}-\lambda_t^{\overline{y},x}\right\vert=O(1)$.  
Hence condition $b_0<1$ is implied by (\ref{chaud}). Theorem \ref{mainstat} leads to the result.$\square$

\subsection{Multinomial logistic autoregressions}
We now provide a multinomial extension of the previous model. We consider the case of a state space $E=\{0,\ldots,N-1\}$ for an integer $N\geq 2$. For $i=1,\ldots,N-1$, assume that
\begin{eqnarray}\label{autoregres2}
\P\left(Y_t=i\vert Y_{t-1}^{-},X_t^{-}\right)&=&\frac{\exp\left(\lambda_{i,t}\right)}{1+\sum_{j=1}^{N-1}\exp\left(\lambda_{j,t}\right)},\\
\lambda_t&=&\sum_{j=1}^q B_j \lambda_{t-j}+\sum_{\ell=1}^p A_{\ell}\overline{Y}_{t-i}+\Gamma X_t,\nonumber
\end{eqnarray}
the $B_j'$s and the $A_{\ell}'$s being matrices of size $(N-1)\times (N-1)$, $\Gamma$ a matrix of size $(N-1)\times d$. Moreover, $\overline{Y}_{t-i}$ takes the $kth$ column of the 
identity matrix $I_{N-1}$ if $Y_{t-i}$ takes the value $k$. In what follows, we denote by $\det(B)$ the determinant of a square matrix $B$. We have the following result.

 \begin{prop}\label{Bin2}
Assume that the covariate process $X$ satisfies {\bf S1} and $\E\log_+\vert X_1\vert$. Assume further that the roots of the polynomial
$$\mathcal{P}(z)=\det\left(I_{N-1}-\sum_{j=1}^q B_j z^j\right)$$
are outside the unit disc. There then exists a unique stationary solution to (\ref{autoregres2}). Moreover, 
the process $\left((Y_t,X_t)\right)_{t\in\Z}$ is stationary and ergodic.
\end{prop} 
\paragraph{Note.} This type of multinomial model is considered in \citet{russell}. Let us point out that 
as for the multinomial regression, a modality of reference is chosen, here $0$. 
In practice, the choice of this modality is often arbitrary and it is an undesirable property 
to have a model depending on this choice. Non-sensitivity to this choice requires that the differences $\lambda_{i,t}-\lambda_{j,t}$ for $i\neq j$ can be obtained via a change in the parameters of the specification of $\lambda_{i,t}$.
This is the case if $B_j=\beta_jI_{N-1}$ in (\ref{autoregres2}), a condition also leading to a more parsimonious model. As for the binary case, more complex models can be obtained by including some interactions between past values of the response and the covariates. Once again, to check our assumptions, it is in general necessary to assume boundedness of this interaction and then boundedness of the process $\lambda_t$.

\paragraph{Proof of Proposition \ref{Bin2}}
Let $F=(F_0,\ldots,F_{N-1}):\R^{N-1}\rightarrow [0,1]^N$ be defined by $F_i(z)=\frac{\exp(z_i)}{1+\sum_{j=1}^{N-1}\exp(z_j)}$ if $1\leq i\leq N-1$ and $F_0(z)=\left(1+\sum_{j=1}^{N-1}\exp(z_j)\right)^{-1}$. Each $F_i$ is a Lipschitz function positive everywhere.
To check Assumption {\bf S2} of Theorem \ref{mainstat}, one can use Lemma \ref{observation-driven} and proceed as for the proof of Proposition \ref{Bin1}.
If $\lambda_t=\left(\mu_t',\ldots,\mu_{t-q+1}'\right)'$, we have $\lambda_t=A\lambda_{t-1}+b_t$
with
$$A=\begin{pmatrix}B_1&\ldots& B_q\\&I_{(q-1)(N-1)}& 0_{(q-1)(N-1),N-1}\end{pmatrix},\quad 
b_t=\begin{pmatrix} \sum_{k=1}^pA_k Y_{t-k}+\Gamma X_t\\0_{(q-1)(N-1),N-1}\end{pmatrix}.$$
The assumption on $\mathcal{P}$ guarantees that the spectral radius of $A$ is less than one (such property is widely known for VAR time series, see for instance \citet{Lut}).
Hence Lemma \ref{observation-driven} guarantees that $Y$ can be seen as chain with complete connections and a geometrically decreasing sequence $(b_m)_m$.  
The crucial point is to check the condition $b_0<1$. Proceeding as in the proof
of Proposition \ref{Bin1}, it is simply necessary to show that for any $c>0$,
\begin{equation}\label{bull} 
\sup_{z,y\in \R^{N-1},\vert y\vert\leq c}\frac{1}{2}\sum_{i=0}^{N-1}\left\vert F_i(z+c)-F_i(z)\right\vert<1.
\end{equation}
Using the equality $\frac{1}{2}\sum_{i=0}^{N-1}\left\vert F_i(z+c)-F_i(z)\right\vert=1-\sum_{i=0}^{N-1}F_i(z+c)\wedge F_i(z)$,
it is enough to show that
\begin{equation}\label{lb}
\alpha:=\inf_{z\in \R^{N-1},\vert y\vert\leq c}\max_{0\leq i\leq N-1}F_i(y+z)\wedge F_i(z)>0.
\end{equation}
Let $k(z)\in \{1,\ldots,N-1\}$ s.t. $z_{k(z)}\geq z_j$ for all $j\in \{1,\ldots,N-1\}$.
We have 
$$F_{k(z)}(y+z)\geq \frac{\exp(z_{k(z)}-c)}{1+(N-1)\exp(z_{k(z)}+c)}\stackrel{k(z)\rightarrow \infty}{\rightarrow} \frac{\exp(-c)}{1+(N-1)\exp(c)},$$
and then 
$$\inf_{z:k(z)\geq 0, \vert y\vert\leq c}\max_{0\leq i\leq N-1}F_i(y+z)\wedge F_i(z)>0.$$
Now, if $k(z)\leq 0$, we have
$$F_0(y+z)\geq \left(1+(N-1)\exp(K)\right)^{-1}$$
and then
$$\inf_{z:k(z)\leq 0, \vert y\vert\leq c}\max_{0\leq i\leq N-1}F_i(y+z)\wedge F_i(z)>0.$$
This shows (\ref{lb}) and then (\ref{bull}) and $b_0<1$. The result of the lemma follows from Theorem \ref{mainstat}.$\square$

\subsection{Discrete choice models}
Here, we assume that $E=\{1,\ldots,N\}$. We want to consider stationary solutions of
\begin{equation}\label{dgraph}
Y_t=\left(\mathds{1}_{\mu_{i,t}+\varepsilon_{i,t}>0}\right)_{1\leq i\leq N},\quad \mu_t=\sum_{j=1}^qB_j \mu_{t-j}+\sum_{k=1}^pA_j Y_{t-j}+\Gamma X_t,
\end{equation}
where $\Gamma$ is a matrix of size $N\times d$ and $A_1,\ldots,A_p,B_1,\ldots,B_q$ are square matrices of size $N\times N$.
Such model is proposed for instance in \citet{Eichler}, \citet{Candelson} or \citet{Nyberg}
for application to financial crisis, business cycles or recession dynamics.
We will show it is possible to construct stationary paths for the dynamic  (\ref{dgraph}) when the two process $X$ and $\varepsilon$ are independent.
More precisely, setting for some $c\in \R^N$ and $I\subset E$, 
$$C_I(c)=\cap_{i\in I}\{\varepsilon_{i,0}>-c_i\}\cap\cap_{i\in E\setminus I}\{\varepsilon_{i,0}\leq -c_i\},$$
we consider stationary processes $(Y_t)_{t\in\Z}$ solution of 
\begin{equation}\label{dgraph2}
\P\left(Y_t=\mathds{1}_I\vert Y_{t-1}^{-},X\right)=\mu\left(C_I(\mu_t)\right),\quad \mu=\P_{\varepsilon_0},
\end{equation}
where $\mathds{1}_I$ is a vector of $\R^N$ with a coordinate $i$ equal to $1$ if $i\in I$ and $0$ otherwise.
\begin{prop}\label{graph}
Assume that $X$ satisfies {\bf A1} with $\E\log_+\vert X_0\vert<\infty$, $\varepsilon_0$ have a distribution with a full support $\R^N$ and a Lipschitz c.d.f. Assume further
that the roots of the polynomial $\mathcal{P}(z)=\det\left(I_N-\sum_{j=1}^qB_jz^j\right)$
are outside the unit disc. There then exists a unique stationary solution to the recursive equations (\ref{dgraph}), (\ref{dgraph2}). Moreover the process $\left((Y_t,X_t)\right)_{t\in\Z}$ is stationary and ergodic. 
\end{prop}

\paragraph{Proof of Proposition \ref{graph}}
Let $I$ be subset of $E$. We will denote by $I^{c}$ the set $E\setminus I$.
As for the previous examples, one can use Lemma \ref{observation-driven} and apply Theorem \ref{mainstat}. The single tricky point is to get condition $b_0<1$.
We only need to show that for any $c\in\R_+^N$,
\begin{equation}\label{lasthope}
\inf_{y,z\in\R^N,\vert y_i\vert\leq c_i,1\leq i\leq N}\sum_{I\subset E}\mu\left(C_I(y+z)\right)\wedge\mu\left(C_I(z)\right)>0.
\end{equation}
If $z\in\R^N$ is fixed, set $I_z=\{i\in E: z_i>0\}$. It is not difficult to show that 
$$\sum_{I\subset E}\mu\left(C_I(y+z)\right)\wedge\mu\left(C_I(z)\right)\geq \mu\left(C_{I_z}\left(-c\mathds{1}_{I_z}+c\mathds{1}_{I^{c}_z}\right)\right).$$
Due to the assumption of full support for $\varepsilon_0$, we have 
$$\mu\left(C_{I_z}\left(-c\mathds{1}_{I_z}+c\mathds{1}_{I_z^c}\right)\right)\leq \inf_{I\subset E}\mu\left(C_I\left(-c\mathds{1}_I+c\mathds{1}_{I^c}\right)\right)>0$$
and (\ref{lasthope}) follows. The result is a then a consequence of Theorem \ref{mainstat}.$\square$

\section{Measures of stochastic dependence}\label{5}

We will now study some weak dependence properties for some processes defined in the previous section.
Many dependence coefficients have been introduced in the literature. See \citet{DandCo} for a survey. 
The notion of strong mixing is probably one of most used for statistical applications. \citet{Doukhan(1994)} is a classical reference on this topic. However, strong mixing conditions are not always easy to check for a bivariate process of type $V_t=(Y_t,X_t)$.  
Under our assumptions, it is possible to show that the conditional probabilities $Y\vert X=x$ satisfy 
$\phi-$mixing conditions. See \citet{FT1} for a discussion in the homogeneous case, the arguments are the same here. However there is no straightforward link between conditional and unconditional mixing. The notion of conditional mixing is considered for instance in \citet{PR}. \citet{CC} give some counterexamples showing that conditional mixing properties do not necessarily entail unconditional strong mixing properties. Moreover, there exist autoregressive processes $(X_t)_{t\in\Z}$ that do not satisfy any strong mixing conditions and alternative dependence coefficients have been proposed in the literature such as the functional dependence of \citet{Wu}, adapted to Bernoulli shifts or the $\tau-$dependence coefficients introduced by \citet{DedP1} and generalized in \citet{DedP2}. 
The latter dependence condition can be used as an alternative to the usual mixing conditions, since the usual deviations inequalities and invariance principles are available for partial sums of $\tau-$dependent sequences. See for instance \citet{DedP1} and \citet{Mer}.  
In this section, we will use our results for bounding either the coefficients of absolute regularity or $\tau-$dependence, depending on the assumption made on the covariate process. 
Since we already pointed out the difficulty of getting unconditional dependence properties from marginal ones, we will assume existence of a particular coupling of the covariate process instead of a particular weak dependence condition. 

\subsection{Absolute regularity and $\tau-$dependence coefficients}

For a stationary process $(V_t)_{t\in\Z}$ taking values in $E\times \R^d$ and $n\in\N^{*}$, we set
$$\beta_V(n)=\E\left[\sup_A\left\vert \P\left((V_n,V_{n+1},\ldots)\in A\vert \mathcal{F}_0\right)-\P\left((V_n,V_{n+1},\ldots)\in A\right)\right\vert\right],$$
where $\mathcal{F}_0=\sigma\left(V_i:i\leq 0\right)$.
We say that $(V_t)_{t\in\Z}$ is absolutely regular or $\beta-$mixing if $\lim_{n\rightarrow \infty}\beta_V(n)=0$.

Next we remind the definition of the coefficients of $\tau-$dependence. 
On $\overline{E}=E\times \R^d$, we consider the distance $\overline{\gamma}$ defined by 
$$\overline{\gamma}(v,v')=\mathds{1}_{v_1\neq v_1'}+\vert v_2'-v_2\vert.$$ 
We also define the following set of Lipschitz functions:
$$\mathcal{L}_{\ell}=\left\{f:\overline{E}^{\ell}\rightarrow \R \mbox{ s.t. } Lip(f):=\sup_{w\neq w'\in \overline{E}^{\ell}}\frac{\left\vert f(w_1,\ldots,w_{\ell})-f(w'_1,\ldots,w'_{\ell})\right\vert}{\sum_{i=1}^{\ell}\overline{\gamma}(w_i,w'_i)}<\infty\right\}.$$
Finally, for a point $v_0\in \overline{E}$, we set 
$$\mathcal{P}_{0,\ell}=\left\{\mu \mbox{ probability measure on } \overline{E}^{\ell}: \int \sum_{j=1}^{\ell}\overline{\gamma}(v_j,v_0)\mu(dv_1,\ldots,dv_{\ell})<\infty\right\}.$$
Note that the set $\mathcal{P}_{0,\ell}$ does not depend on the point $v_0$. Next, for $\mu,\nu\in \mathcal{P}_{0,\ell}$, we define
$$\mathcal{W}_{1,\ell}(\mu,\nu)=\sup_{f\in\mathcal{L}_{\ell}}\left\{\int f d\mu-\int fd\nu: Lip(f)\leq 1\right\}.$$
Remember that from Kantorovich's duality, we have 
$$\mathcal{W}_{1,\ell}(\mu,\nu)=\inf\left\{\int \sum_{j=1}^{\ell}\overline{\gamma}(v_j,v'_j)\Gamma(dv_1,\ldots,dv_{\ell},dv_1',\ldots,dv_{\ell}')\right\}$$
where the infinimum is taken on the set of probability measures $\Gamma$ on $\overline{E}^{\ell}\times \overline{E}^{\ell}$ having marginals $\mu$ and $\nu$. 

For $t\in\Z$, let $V_t=(Y_t,X_t)$, $\mathcal{F}_t=\sigma\left((\varepsilon_j,Y_j): j\leq t\right)$.
For an integer $\ell\geq 1$ and $j_1<\cdots<j_{\ell}$ in $\Z$, set $J=\{j_1,\ldots,j_{\ell}\}$ and 
$$U_J=\left(V_{j_1},V_{j_2},\ldots,V_{j_{\ell}}\right).$$
According to \citet{DedP1}, we define $\tau-$dependence coefficients between $U_J$ and $\mathcal{F}_0$ by 
$$\tau\left(\mathcal{F}_0,U_J\right)=\E\mathcal{W}_{1,\ell}\left(\P_{U_J\vert\mathcal{F_0}},\P_{U_J}\right),$$
where $\P_{U_J\vert\mathcal{F}_0}$ denotes the conditional distribution of $U_J$ given $\mathcal{F}_0$. 
We then define for some integers $k\geq 1$ and $n\geq 1$,
$$\tau^{(k)}(n)=\max_{1\leq \ell\leq k}\frac{1}{\ell}\sup\left\{\tau\left(\mathcal{F}_0,U_{j_1,\ldots,j_{\ell}}\right), n\leq j_1<\cdots<j_{\ell}\right\},$$
and $\tau_V(n)=\sup_{k\geq 1}\tau^{(k)}(n)$. 

Note that the initial definition of the $\tau-$dependence coefficients defined in \citet{DedP1} were  defined when the distance $\overline{\gamma}$ is the $\ell_1-$metric. However, as $E$ is a finite set,
the two metric are equivalent. Indeed, one can always code the elements of the finite set $E$ as vectors of the canonical basis of $\R^N$ ($N$ is the number of elements of $E$) and in this case, we have, for $x,y\in E$, 
$$\frac{1}{2}\sum_{i=1}^N\vert x_i-y_i\vert=\mathds{1}_{x\neq y}.$$
Then one can assume that the process $(V_t)_{t\in \Z}$ takes values in $\R^{N+d}$ and choose $\overline{\gamma}$ as the corresponding $\ell^1-$metric. In this case, \citet{DedP1} developed various limiting theorem for partial sums, when $\tau_V(n)\rightarrow 0$.   
Let us also mention that such dependence coefficients were generalized in \citet{DedP2}, when the state space $(\overline{E},\overline{\gamma})$ is a general Polish space.

\subsection{Control of dependence coefficients}

For bounding the dependence coefficients defined in the previous section, 
we will assume a representation of the form
$X_t=g\left(S_t\right)$, where $(S_t)_{t\in\Z}$ is a stationary Markov chain taking values in a Polish space $\mathcal{S}$ and $g:\mathcal{S}\rightarrow \R^d$ is a measurable function.
We will denote by $\gamma$ a metric on $\R^d$ that will be either the discrete metric, i.e. $\gamma(x_1,x_2)=\mathds{1}_{x_1\neq x_2}$, or the $\ell_1-$metric, i.e. $\gamma(x_1,x_2)=\vert x_2-x_1\vert:=\sum_{i=1}^d\vert x_{i,1}-x_{i,2}\vert$ for $x_j=(x_{1,j},\ldots,x_{d,j})\in\R^d$, $j=1,2$.  

We introduce the probability kernel $P$ from $\left(\mathcal{S},\mathcal{B}(\mathcal{S})\right)$
to $\left(\mathcal{S}^{\N^*},\mathcal{B}\left(\mathcal{S}^{\N^*}\right)\right)$
and such for $s_0\in\mathcal{S}$, $P(s_0,\cdot)$ is the distribution of $(S_t)_{t\geq 1}$ conditional on $S_0=s_0$. Let also $\pi$ be the invariant probability of the chain.
We set $\Omega^{*}=\mathcal{S}^{\N^*}\times \mathcal{S}^{\N^*}$. For $t\in\N^*$, we denote by $S_{1,t}:\Omega\rightarrow \mathcal{S}$ and $S_{2,t}:\Omega^*\rightarrow \mathcal{S}$ the coordinate applications $S_{1,t}((s,\overline{s}))=s_t$ and $S_{2,t}((s,\overline{s}))=\overline{s}_t$.

\begin{description}
\item [S1']  We assume that 
there exists a probability kernel $\widetilde{P}$ from $\left(\mathcal{S}^2,\mathcal{B}\left(\mathcal{S}^2\right)\right)$ to $\left(\Omega^{*},\mathcal{B}\left(\Omega^*\right)\right)$ such that $\widetilde{P}\left((s_0,\overline{s}_0),\cdot\right)$ is a coupling of $P(s_0,\cdot)$ and $P(\overline{s}_0,\cdot)$  and 
$$a_t=:=\int d\pi(s_0)d\pi(\overline{s}_0)\widetilde{E}_{s_0,\overline{s_0}}\left[\gamma\left(g(S_{1,t}),g(S_{2,t})\right)\right]$$
satisfies $\lim_{t\rightarrow \infty}a_t=0$. $\widetilde{E}_{s_0,\overline{s_0}}$ denotes the expectation under $\widetilde{P}\left((s_0,\overline{s}_0),\cdot\right)$.

\item[S3] There exists a sequence $(e_m)_{m\in \N}$ such that $\sum_{m\geq 1}e_m<\infty$ and for all $(y,z,z')\in G^{\N}\times E^{\N}\times E^{\N}$, 
$$d_{TV}\left(q\left(\cdot\vert y,z\right),q\left(\cdot\vert y,z'\right)\right)\leq \sum_{i\geq 0}e_i\vert z_i-z_i'\vert.$$  
\end{description}

\paragraph{Notes}
\begin{enumerate}
\item
Assume that $\gamma$ is induced by the $\ell_1-$norm on $\R^d$ and the covariate process $X$ is a Bernoulli shift. We remind that a process $X$ is called a 
a Bernoulli shift if there exists a measurable space $\Lambda$, a random sequence $\varepsilon\in \Lambda^{\Z}$ 
of i.i.d. random variables and a measurable application $g:\Lambda^{\N}\rightarrow G$ such that 
$$X_t=g\left(\varepsilon_t,\varepsilon_{t-1},\ldots\right),\quad t\in \Z.$$
We point out that such a representation is valid for many time series models found in the literature from linear processes of ARMA type to GARCH processes.
In this case we set $S_t=\left(\varepsilon_t,\varepsilon_{t-1},\ldots,\right)$ which takes values in $\mathcal{S}=\Lambda^{\N}$. Here, the kernel $\widetilde{P}$ in Assumption {\bf S1'} can be defined as the probability distribution of $\left(\varepsilon^{(s_0)},\varepsilon^{(\overline{s}_0)}\right)$
where for $s\in \Lambda^{\N}$, $\varepsilon_t^{(s)}=\varepsilon_t$ for $t\geq 1$ and $\varepsilon_t^{(s)}=s_{-t}$ for $t\leq 0$. In this case, we have the expression $a_t=\E\left(\left\vert X_t-\overline{X}_t\right\vert\right)$ where $t\geq 1$, 
$$\overline{X}_t=g\left(\varepsilon_t,\ldots,\varepsilon_1,\varepsilon'_0,\varepsilon'_{-1},\ldots\right)$$
and $\varepsilon'$ is an independent copy of $\varepsilon$. A martingale argument shows that the condition $\lim_{t\rightarrow \infty}a_t=0$ is automatically satisfied.

\item
When $\gamma(x_1,x_2)=\mathds{1}_{x_1\neq x_2}$,  Assumption {\bf S1'} is implied by the existence of a so-called successful coupling of two chains with different initial values. See \citet{Lind} for a discussion of existence of successful coupling for Markov chains and in particular Theorem $14.10$ which shows an equivalence with the weak ergodicity property of the Markov chain.
Existence of a successful coupling means that there exists a probability kernel $\widetilde{P}$, defined as in Assumption {\bf S1'} and such that 
$$\widetilde{P}\left((s_0,\overline{s}_0),\left\{\exists n_0\in\N: S_{1,n}=S_{2,n}, n\geq n_0\right\}\right)=1.$$ Then the condition $\lim_{t\rightarrow \infty}a_t=0$ is automatically satisfied for such kernel. Indeed, we have 
$$\widetilde{E}_{s_0,\overline{s}_0}\left[\gamma\left(g(S_{1,t},g(S_{2,t})\right)\right]\leq \widetilde{P}\left((s_0,\overline{s}_0),\cup_{i\geq t}\{S_{1,i}\neq S_{2,i}\}\right)\rightarrow 0.$$
and we get $a_t\rightarrow 0$ from Lebesgue's theorem.  

\item
One can also assume that $S_t=\left(\varepsilon_t,\varepsilon_{t-1},\ldots\right)$
where $\left(\varepsilon_t\right)_{t\in \Z}$ is a homogenous chain with complete connections satisfying Assumption {\bf S2}. Using Lemma \ref{central}, one can check Assumption {\bf A1'} when $g(S_t)=h(\varepsilon_t,\ldots,\varepsilon_{t-k})$ for some integer $k$ and function $h$.
\end{enumerate}

In what follows, we denote by $\L^p$ the Lebesgue space of random variables taking values in $\R^d$ and possessing a moment of order $p$ (or bounded a.s. if $p=\infty$) and $\Vert\cdot\Vert_p$ the corresponding norm. 
\begin{theo}\label{mixing1}
Assume that Assumptions {\bf S1-S1'} and {\bf S2-S3} hold true. 
\begin{enumerate}
\item
Assume that $\gamma$ is the discrete metric and that $X_0\in \L^p$ for some $p\in [1,\infty]$. Set $q=\frac{p}{p-1}$ and $c_t=\max(1,2\Vert X_0\Vert_p)a_t^{1/q}$. Then, for $n\in\N^{*}$, we have 
$$\beta_V(n)\leq \sum_{j\geq n}g_j,$$
with $g_j=b_{j-1}^{*}+c_j+\kappa_j+\sum_{i=0}^{j-2}b_i^{*}\kappa_{j-i-1}$ and $
\kappa_j=\sum_{s=0}^{j-1} e_s c_{j-s}+2\sum_{s\geq j} e_s\E\left[\vert X_0\vert\right]$.
\item
Assume that $\gamma$ is induced by the $\ell_1-$norm. Then, for $n\in\N^{*}$, we have 
$$\tau_V(n)\leq \sup_{j\geq n}h_j,$$
with
$h_j=b_{j-1}^{*}+a_j+\kappa_j+\sum_{i=0}^{j-2}b_i^{*}\kappa_{j-i-1}$, 
$\kappa_j=\sum_{s=0}^{j-1} e_s a_{j-s}+2\sum_{s\geq j} e_s\E\left[\vert X_0\vert\right]$.
\end{enumerate}
\end{theo}

\paragraph{Notes}
\begin{enumerate}
\item
Note that under our assumptions, $\lim_{i\rightarrow \infty}g_i=0$. This is essentially due to the fact that if $(u_n)_{n\geq 0}$ is a summable sequence of nonnegative real numbers and $(v_n)_{n\in\N}$ a sequence of nonnegative real numbers
converging to $0$, then $\lim_{n\rightarrow \infty}\sum_{i=0}^nu_iv_{n-i}=0$.
Summability....

\item
When $b_m=O(m^{-k})$ for some $k\in\N^*$, then we also have $b_m^{*}=O(m^{-k})$. See Lemma \ref{moment} given in the Appendix. 
\end{enumerate}

\paragraph{Proof of Theorem \ref{mixing1}}
Using Theorem \ref{mainstat}, one can define the process $Y$ conditional on $S$, instead of conditional on $X$. The resulting process will be the unique stochastic process satisfies (\ref{spec}). 
To this end, we simply change the set $\mathcal{D}$ by the set $\mathcal{D}_0=\left\{s\in\mathcal{S}^{\Z}:
\left(g(s_t)\right)_{t\in\Z}\in \mathcal{D}\right\}$. One can then consider that the distribution of the Markov chain $S$ is supported on $\mathcal{D}_0$.
Secondly, for bounding $\beta_V(n)$ or $\tau_V(n)$, one can replace the sigma-field $\mathcal{F}_0$ by a larger one. This follows from the properties of the conditional expectations.  
Let $\mu$ be the probability distribution of $S_0^{-}$. We also denote by $K_{s_0^{-}}$ the probability distribution of 
$Y_0^{-}$ conditional on $S=s$. Note that from Theorem \ref{mainstat}, this conditional distribution only depends on $s_0^{-}$. 

For $(t,w,z,s,\overline{s})\in \Z\times E\times E^{\N}\times \mathcal{D}_0\times\mathcal{D}_0$ and $x=x(s,\overline{s})=\left((g(s_t))_{t\in\Z},(g(s_t))_{t\in\Z}\right)$, we set
$$q_t^{x}(w\vert z)=q\left(w\vert z,g(s_t)^{-}\right),\quad \overline{q}_t^{x}(w\vert z)=q\left(w\vert z,g(\overline{s}_t)^{-}\right).$$
Next setting $s_1^{+}=(s_1,s_2,\ldots)$ and $y_1^{+}=(y_1,y_2,\ldots)$ for any $(s,y)\in \mathcal{D}_0\times E^{\N}$, 
we consider a probability measure $\P$ on $\Omega=\mathcal{S}^{\Z}\times\mathcal{S}^{\Z}\times E^{\Z}\times E^{\Z}$ endowed with its Borel $\sigma-$field and defined by
\begin{eqnarray*}
&&\P\left(ds,d\overline{s},dy,d\overline{y}\right)\\
&=&\mu\left(ds_0^{-}\right)\mu\left(d\overline{s}_0^{-}\right)K_{s_0^{-}}\left(dy_0^{-}\right)K_{\overline{s}_0^{-}}\left(d\overline{y}_0^{-}\right)
\widetilde{P}\left((s_0,\overline{s}_0), (ds_1^{+},d\overline{s}_1^{+})\right)
\widetilde{Q}_{0,x(s,\overline{s}),z,\overline{z}}\left(dy_1^{+},d\overline{y}_1^{+}\right).
\end{eqnarray*}
We remind that $\widetilde{P}$ is defined in Assumption {\bf S1'} and $\widetilde{Q}$ is defined in Lemma \ref{central}.
On $\Omega$, we will still denote, for $t\in \Z$,  the coordinate applications by $Y_t,\overline{Y}_t,S_t,\overline{S}_t$.
Let us also point out that the measure 
$$\widetilde{P}\left((s_0,\overline{s}_0), (ds_1^{+},d\overline{s}_1^{+})\right)
\widetilde{Q}_{0,x(s,\overline{s}),z,\overline{z}}\left(dy_1^{+},d\overline{y}_1^{+}\right)$$
is a coupling of two conditional distributions, the distribution of $\left((Y_t,S_t)\right)_{t\geq 1}$ conditional on $Y_j=y_j,S_j=s_j$ for $j\leq 0$ and 
the distribution of $\left((Y_t,S_t)\right)_{t\geq 1}$ conditional on $Y_j=\overline{y}_j,S_j=\overline{s}_j$ for $j\leq 0$.   
Next we set $\mathcal{G}_0=\sigma\left((Y_j,S_j): j\leq 0\right)$ and $\overline{\mathcal{G}}_0=\sigma\left((\overline{Y}_j,\overline{S}_j): j\leq 0\right)$.
Let $J=\{j_1,\ldots,j_{\ell}\}\subset \N^{*}$.
Note that the two sigma-fields $\mathcal{G}_0$ and $\overline{\mathcal{G}}_0$ are independent.

\begin{enumerate}
\item
For the absolute regularity coefficients, 
we use the bounds
\begin{eqnarray*}
\beta_V(n)&\leq& \E\left[\sup_A\left\vert \P\left((V_n,V_{n+1},\ldots)\in A\vert \mathcal{G}_0\right)-
\P\left((\overline{V}_n,\overline{V}_{n+1},\ldots)\in A\vert \overline{\mathcal{G}}_0\right)\right\vert\right]\\
&\leq& \sum_{t\geq n}\E\left[\P\left(V_n\neq \overline{V}_n\vert \mathcal{G}_0\vee \overline{\mathcal{G}}_0\right)\right].
\end{eqnarray*}
Next, using Lemma \ref{central} and Assumption {\bf S3}, we have
\begin{eqnarray*}
\P\left(Y_t\neq \overline{Y}_t\vert \mathcal{G}_0\vee\overline{\mathcal{G}}_0\vee\sigma(S_1^{+},\overline{S}_1^{+})\right)&\leq& \sup_{z_0^{-},\overline{z}_0^{-}}\widetilde{Q}_{0,x(S,\overline{S}),z_0^{-},\overline{z}_0^{-}}\left(\{y_t\neq \overline{y}_t\}\right)\\
&\leq& b_{t-1}^{*}+ \sup_{g\in E^{\N}}d_{TV}\left(q(\cdot\vert g,S_t^{-}),q(\cdot\vert g,\overline{S}_t^{-})\right)\\
&+&\sum_{\ell=0}^{t-2}b_{\ell}^{*}\sup_{g\in E^{\N}}d_{TV}\left(q(\cdot\vert g,S_{t-\ell-1}^{-}),q(\cdot\vert g,\overline{S}_{t-\ell-1}^{-})\right)\\
&\leq& b_{t-1}^{*}+\sum_{i\geq 0}e_i G_{t-i}+\sum_{\ell=0}^{t-2}b_{\ell}^{*}\sum_{i=0}^{\infty}e_i G_{t-\ell-i-1},
\end{eqnarray*}
where for any $t\in \Z$, $G_t=\left\vert g(S_t)-g(\overline{S}_t)\right\vert$.
Using Holder inequality, we have
$$\E(G_t)=\E\left[\left\vert g(S_t)-g(\overline{S}_t)\right\vert \mathds{1}_{g(S_t)\neq g(\overline{S}_t)}\right]\leq 2\Vert X_0\Vert_p \P\left(g(S_t)\neq g(\overline{S}_t)\right)^{1/q}.$$
From the definition of the coupling, we have $\P\left(g(S_t)\neq g(\overline{S}_t)\right)=a_t$ where $a_t$ is defined in {\bf S1'}. Since,
$$\P\left(V_n\neq \overline{V}_n\vert \mathcal{G}_0\vee \overline{\mathcal{G}}_0\right)\leq
\P\left(Y_t\neq \overline{Y}_t\vert \mathcal{G}_0\vee\overline{\mathcal{G}}_0\right)+\P\left(g(S_t)\neq g(\overline{S}_t\vert \sigma(S_0,\overline{S}_0)\right)$$
and $a_t=\P\left(g(S_t)\neq g(\overline{S}_t)\right)\leq c_t$,
the bound for $\beta_V(n)$ follows after integrating the previous inequalities.

\item
We have
$$\tau\left(\mathcal{F}_0,U_J\right)\leq \E\left[W_{1,\ell}\left(\P_{U_J\vert\mathcal{F}_0},\P_{\overline{U}_J\vert \overline{\mathcal{F}}_0}\right)\right],$$
with $\overline{U}_J=\left(\overline{V}_{j_1},\ldots,\overline{V}_{j_{\ell}}\right)$ and $\overline{V}_t=\left(\overline{Y}_t,g\left(\overline{S}_t\right)\right)$ for $t\in\Z$.
Using our coupling we have
\begin{eqnarray*}
&&W_{1,\ell}\left(\P_{U_J\vert\mathcal{G}_0},\P_{\overline{U}_J\vert \overline{\mathcal{G}}_0}\right)\\
&\leq& \sum_{i=1}^{\ell}\left[\P\left(Y_{j_i}\neq \overline{Y}_{j_i}\vert \mathcal{G}_0\vee\overline{\mathcal{G}}_0\right)+\E\left(\vert g(S_{j_i})-g(\overline{S}_{j_i})\vert \big\vert \sigma\left(S_0,\overline{S}_0\right)\right)\right]\\
&\leq& \ell\sup_{t\geq n}\left[\P\left(Y_t\neq \overline{Y}_t\vert \mathcal{G}_0\vee\overline{\mathcal{G}}_0\right)+\E\left(\vert g(S_t)-g(\overline{S}_t)\vert \big\vert \sigma\left(S_0,\overline{S}_0\right)\right)\right].\\
\end{eqnarray*}
From the definition of our coupling and Assumption {\bf S1'}, we have 
$$\E\left[\vert g(S_t)-g(\overline{S}_t)\right]\leq a_t.$$
Next one can bound $\P\left(Y_t\neq \overline{Y}_t\right)$
as in the previous point and we have directly $\E(G_t)=a_t$.
The proposed upper-bound for $\tau_V(n)$ easily follows.$\square$
\end{enumerate}

Next we give a result focused on observation-driven models. Our conditions will be specified for the examples already mentioned in Section \ref{4}

\begin{cor}\label{import}
Assume that Assumptions {\bf S1-S1'} hold true.
Suppose that $q\left(\cdot\vert Y_{t-1}^{-},X_t^{-}\right)=\check{q}\left(\cdot\vert \lambda_t\right)$ is Lipschitz in $\lambda_t$ with $\lambda_t=G_{Y_{t-1},X_t}\left(\lambda_{t-1}\right)$ and the assumptions of Lemma \ref{observation-driven} are satisfied.
\begin{enumerate}
\item
Assume that {\bf S1'} holds true for the discrete metric.
If $a_i=O(i^{-\kappa})$ with $\kappa>q$, $X_0\in \L_p$ and $p^{-1}+q^{-1}=1$, we have $\beta_V(n)=O\left(n^{-\kappa/q+1}\right)$ and the process $V$ is absolutely regular. If $a_i=O\left(\rho^i\right)$ for some $\rho\in (0,1)$, then $\beta_V(n)=O\left(\overline{\rho}^n\right)$ for some $\overline{\rho}\in (0,1)$.
\item
Assume that {\bf S1'} holds true for the $\ell_1-$metric.
If $a_i=O(i^{-\kappa})$ with $\kappa>1$, we have $\tau_V(n)=O\left(n^{-\kappa}\right)$. If $a_i=O\left(\rho^i\right)$ for some $\rho\in (0,1)$, then $\tau_V(n)=O\left(\overline{\rho}^n\right)$ for some $\overline{\rho}\in (0,1)$.
\end{enumerate}
\end{cor}

\paragraph{Proof of Corollary \ref{import}}
From Lemma \ref{observation-driven}, the coefficients $b_i$ decay exponentially and from Lemma \ref{moment}, so do 
the corresponding coefficients $b_i^{*}$. 
Moreover, the coefficients $e_i$ in Assumption {\bf S3} also decay exponentially fast.
For a polynomial decay, we have $g_j=O\left(j^{-\kappa/q}\right)$ and $h_j=O\left(j^{-\kappa}\right)$ in Theorem \ref{mixing1}. The result of the corollary then follows. $\square$

\subsection{Heredity}
We now provide an upper bound for the $\tau-$dependence coefficients of some functionals of 
the process $(V_t)_{t\in\Z}$ that will be quite general for statistical applications.
More precisely, our aim is to get limit theorems or deviation inequalities for some partial sums of the form
$$S_n=\frac{1}{n}\sum_{t=1}^n W_t,\quad W_t:=f\left(V_t,V_{t-1},\ldots\right),$$
for some suitable functions $f$ and which can include functionals of $(Y_{t-j},\lambda_{t-j})_{j\geq 0}$ as a special case, where $\lambda$ denotes the latent process of observation-driven models satisfying our assumptions.
We point out that, due to the discrete nature of the process $(Y_t)_{t\in\Z}$, the process $(\lambda_t)_{t\in\Z}$ is not necessarily absolutely regular. \citet{Neu} studied this problem for the Poisson autoregressive process. Then one can not deduce absolute regularity of the process $(W_t)_{t\in\Z}$ from that of $(V_t)_{t\in\Z}$. In contrast, such heredity is possible for the $\tau-$dependence. Since we did not find a precise reference for such a property, we state a result in the present paper. 
In this section, we assume that $\E\vert X_0\vert<\infty$.
For a sequence of positive and summable coefficients $(\alpha_i)_{i\geq 1}$, we set
$$H_{\alpha}=\left\{\left((y_i,x_i)\right)_{i\geq 0}\in (E\times\R^d)^{\N}: \sum_{i\geq 0}\alpha_i\vert x_i\vert<\infty\right\}.$$
For $z=(z_i)_{i\in\N}$, we alo set $\vert z\vert_{\alpha}=\sum_{i\geq 0}\alpha_i \vert z_i\vert$. 
Now we consider a process $(W_t)_{t\in \Z}$ defined by  
$$W_t=f\left(V_t,V_{t-1},\ldots\right),\quad t\in\Z,$$
where $f:H_{\alpha}\rightarrow \R^k$ satisfies for $v=((y_i,x_i))_{i\geq 0}$ and $\overline{v}=((\overline{y}_i,\overline{x}_i))_{i\geq 0}$ in $H_{\alpha}$, 
$$\left\vert f(v)-f(\overline{v})\right\vert\leq \left(1+\vert x\vert_{\alpha}^{p-1}+\vert y\vert_{\alpha}^{p-1}\right)\cdot\vert x-y\vert_{\alpha},$$
The $\tau-$dependence coefficients for $W$ are defined as that of $V$, replacing the space $\overline{E}$ by $\R^k$ and the metric $\overline{\gamma}$ on $\overline{E}$ by a norm, still denoted by $\vert\cdot\vert$ on $\R^k$.
We remind that for $t\in\Z$, $\mathcal{F}_t=\sigma\left((Y_i,S_i): i\leq t\right)$.

\begin{prop}\label{heredity}
Assume that $\sum_{i\geq j+1}\alpha_i=O\left(j^{-\eta}\right)$ and $\tau_V(i)=O\left(i^{-\kappa}\right)$ for $\eta,\kappa>1$. Then, if $\E\left[\vert X_0\vert^{p+q}\right]<\infty$ for some $q>0$, we have the bound 
$$\tau_W(i)=O\left(i^{-\kappa'}\right),\quad \kappa'=\min\left(\eta-1,(\kappa-1)\frac{q+2}{q+p+1}\right).$$
\end{prop}

\paragraph{Proof of Proposition \ref{heredity}}
Let $T$ be a positive real number and $j$ a positive integer.
We introduce two other sequences of random variables. For $t\in\Z$, let
$$W_t^{(1)}=f\left(\phi_T(V_t),\phi_T(V_{t-1}),\ldots\right),$$
where for $v=(y,x)\in E\times\R^d$,
$\phi_T(v)=\left(y,\left((-T)\vee x_i\wedge T\right)_{1\leq i\leq d}\right)$. Moreover let
$$W^{(2)}_t=f\left(\phi_T(V_t),\ldots,\phi_T(V_{t-j}),0,0,\ldots\right).$$
Since we use the same filtrations for evaluating $\tau_{W},\tau_{W^{(1)}}$ and $\tau_{W^{(2)}}$, one can use the following inequalities which are a consequence of the definition of the Wasserstein metric.
$$\tau_{W}(i)\leq \tau_{W^{(1)}}(i)+2\E\left\vert W_t-W^{(1)}_t\right\vert\leq \tau_{W^{(2)}}(i)+2\E\left\vert W^{(1)}_t-W^{(2)}_t\right\vert+2\E\left\vert W_t-W^{(1)}_t\right\vert.$$
First, it is easily seen that 
$$\E\left\vert W^{(1)}_t-W^{(2)}_t\right\vert=O\left(\sum_{i\geq j+1}\alpha_i\right).$$
Moreover, the function $(v_1,\ldots,v_j)\mapsto f\left(\phi_T(x_1),\ldots,\phi_T(x_j),0,\ldots\right)$ is Lipschitz with a Lipschitz constant bounded by $T^{p-1}$ (up to a constant). This gives the bound
$$\tau_{W^{(2)}}(i)=O\left(T^{p-1}j\tau_{V}(i-j)\right).$$
Moreover, we have
$$\E\left\vert W_t-W_t^{(1)}\right\vert\leq 2\sum_{i\geq 1}\E\left(1+\vert V_{n,t+1}^{-}\vert_{\alpha}^{p-1}\right)\cdot \vert V_{n,t+1-i}\vert \mathds{1}_{\vert V_{n,t+1-i}\vert\geq T}.$$
From the moment assumption on $X_t$, we get 
$$\E\left\vert W_t-W^{(1)}_t\right\vert=O(T^{-q}).$$
The result then follows by choosing $j=[i/2]$ and $T=\left(j\tau_{V}(i-j)\right)^{-\frac{1}{p+q+1}}$.$\square$

\section{Perspectives in statistics}\label{6}

Let us now give some possible applications of our results to statistical inference in the models presented in Section \ref{4}.
\begin{enumerate}
\item
The first problem concerns parametric estimation which has been extensively studied for other observation-driven models. Usually, only ergodicity is necessary 
to get consistency and asymptotic normality of the conditional likelihood estimator. For instance, \citet{DDM} studied this problem when some general observation-driven models are well specified and misspecified.
Then one could obtain similar results for our models satisfying the assumptions {\bf S1-S2}. 
Let us also mention that models of infinite order can be considered, such as (\ref{autoregres0}), with a parametric form for the parameters $a_j=a_j(\theta)$ and a decay in the dependence which is not exponential (in contrast to observation-driven models). For instance $a_j(\theta)=\theta_1 j^{-\theta_2}$ with $\theta_2>2$. 
\item
The second problem concerns discrete choice models as in (\ref{dgraph}) and the estimation of the distribution of $\varepsilon_0$ modeled via a parametric copula, as in \citet{Eichler}. For instance, the distribution of $\varepsilon_t$ can be marginally Gaussian or logistic and we obtain a multivariate version of the univariate probit or logistic binary model.
Our results can then be used to solve the issues mentioned in their paper. 
\item
The third problem concerns semi-parametric estimation in our models. Recently, \citet{Park} investigated this problem for finite-order models.
Let us describe an approach for the binary time series models given in Section \ref{4}. Our aim is to estimate the function $F$ as well as a vector $\theta$ of autoregressive parameters. One can then 
maximize 
$$\theta\mapsto \sum_{t=1}^n\left[Y_t \log \hat{F}_{\theta}\left(\mu_t(\theta)\right)+(1-Y_t)\log\left(1-\hat{F}_{\theta}\left(\mu_t(\theta)\right)\right)\right],$$
with 
$$\hat{F}_{\theta}(z)=\frac{\sum_{t=1}^nY_tK_h\left(z-\mu_t(\theta)\right)}{\sum_{t=1}^nK_h\left(z-\mu_t(\theta)\right)}$$    
and $K$ is a kernel, $h>0$ a bandwidth parameter and $K_h=h^{-1}K(\cdot/h)$. If $\hat{\theta}$ is such maximizer, $\P\left(Y_t=1\vert Y_{t-1}^{-},X_t^{-}\right)$ can be estimated by $\hat{F}_{\hat{\theta}}\left(\mu_t(\hat{\theta})\right)$. The dependence properties stated in Section \ref{5} will be essential to derive asymptotic properties of this estimator. With respect to the problem considered in \citet{Park}, the main interest of this semi-parametric approach is that one can use much more lags values for the response and the covariates in order to predict the $Y_t$'s.
\end{enumerate}

\section{Appendix}
\subsection{Auxiliary lemmas}

\begin{lem}\label{auxilliaire}
Suppose that $\left\{q_1,q_2,\ldots,\right\}$ is a family of kernels from $\left(E^{\N},\mathcal{B}(E^{\N})\right)$
to $\left(E,\mathcal{B}(E)\right)$ satisfying the Assumption {\bf A}. For $x,y\in E^{\N}$, 
let $\left(U_n^{x,y},V_n^{x,y}\right)_{n\in \Z}$ be the coupling such that  
$$U_n^{x,y}=x_{-n},\quad V_n^{x,y}=y_{-n},\quad n\geq 0$$
and for $n\geq 1$,
$$\P\left(U_n^{x,y}\neq V_n^{x,y}\vert U_{n-j}^{x,y},V_{n-j}^{x,y}; j\geq 1\right)=d_{TV}\left[q_n\left(\cdot\vert U^{x,y}_{n-j};j\geq 1\right),q_n\left(\cdot\vert V^{x,y}_{n-j}; j\geq 1\right)\right].$$   
Then we have the following bounds.
\begin{enumerate}
\item
For $n\geq 1$, $\P\left(U_n^{x,y}\neq V_n^{x,y}\right)\leq b_{n-1}^{*}$,
where $b^{*}$ is defined in the statement of Lemma \ref{central}.
\item
For $n,k\geq 1$, 
$$\P\left(\left(U_n^{x,y},\ldots,U_{n+k}^{x,y}\right)\neq \left(V_n^{x,y},\ldots,V_{n+k}^{x,y}\right)\right)\leq \sum_{j=n-1}^{n+k-1}b_j^{*}.$$
\end{enumerate}
\end{lem}

\paragraph{Proof of Lemma \ref{auxilliaire}}
\begin{enumerate}
\item
For the case of homogeneous transitions, i.e. $q_t$ does not depend on $t$, this result is proved in \citet{Bres}, Proposition $1$.
The nonhomogeneous case follows similarly but for sake of completeness, we give below a complete proof. 
First, we note that if $x_0^{-}=y_0^{-}$, Assumption {\bf A} guarantees that $\P\left(U_n^{x,y}\neq V_n^{x,y}\right)=0$ for all $n\geq 1$.
Next, assume that $x_0^{-}\neq y_0^{-}$. We now omit the superscript $x,y$ because all the bounds will be uniform with respect to $x$ and $y$.  
For a positive integer $n$, we set
$$T_n=\inf\left\{m\geq 0: U_{n-m}\neq V_{n-m}\right\}.$$
Note that $T_n$ is finite a.s. Moreover, $\P\left(U_n\neq V_n\right)=\P\left(T_n=0\right)$.
One can also note that $\P(T_1=0)\leq b_0=b_0^{*}$ by the definition of the coupling. 
To show that $\P\left(T_n=0\right)\leq b_{n-1}^{*}$ for $n\geq 2$, we will show that
\begin{equation}\label{goodbound}
\P\left(S_{n-1}^{(b)}\geq k\right)\leq \P\left(T_n\geq k\right),\quad n\geq 1, k\geq 0.
\end{equation}
This will follow from an induction argument.
For $n=1$, the result follows from the equality $S_0^{(b)}=0$. 
Let us now assume that (\ref{goodbound}) is true for some $n\geq 1$. 
The inequality is automatic for $k=0$, so let us assume that $k\geq 1$. 
If $j\geq 1$, we have the inclusion $\{T_{n+1}=j\}\subset\{T_n=j-1\}$. 
Moreover, on the event $\{T_n=j-1\}$, the possible values for $T_{n+1}$
are $j$ or $0$. From the definition of maximal coupling and the $b_n'$s, we have 
$$\P\left(T_{n+1}=0,T_n=j-1\right)=\E\left[\P\left(U_{n+1}\neq V_{n+1}\vert W_n^{-}\right)\mathds{1}_{T_n=j-1}\right]\leq b_{j-1}\P\left(T_n=j-1\right).$$
We then deduce that 
$$\P\left(T_{n+1}=j\right)=\P\left(T_{n+1}=j,T_n=j-1\right)\geq (1-b_{j-1})\P\left(T_n=j-1\right).$$
Now, we use the bound
$$\P\left(T_{n+1}\geq k\right)=\sum_{j\geq k}\P\left(T_{n+1}=j\right)\geq \sum_{j\geq k}(1-b_{j-1})\P\left(T_n=j-1\right).$$
One can show that
$$\sum_{j\geq k-1}(1-b_j)\P\left(T_n=j\right)=(1-b_{k-1})\P\left(T_n\geq k-1\right)+\sum_{j\geq k}\left(b_{j-1}-b_j\right)\P\left(T_n\geq j\right).$$
We then use the induction hypothesis and the fact the sequence $(b_n)_{n\geq 0}$ is nonincreasing to get 
$$\P\left(T_{n+1}\geq k\right)\geq (1-b_{k-1})\P\left(S^{(b)}_n\geq k-1\right)+\sum_{j\geq k}\left(b_{j-1}-b_j\right)\P\left(S^{(b)}_n\geq j\right).$$
Using similar computations, we have 
$$(1-b_{k-1})\P\left(S^{(b)}_n\geq k-1\right)+\sum_{j\geq k}\left(b_{j-1}-b_j\right)\P\left(S^{(b)}_n\geq j\right)=\sum_{j\geq k-1}(1-b_j)\P\left(S_{n-1}^{(b)}=j\right).$$
Moreover, 
\begin{eqnarray*}
\sum_{j\geq k-1}(1-b_j)\P\left(S_{n-1}^{(b)}=j\right)&=&\sum_{j\geq k-1}\P\left(S^{(b)}_{n-1}=j\right)P(j,j+1)\\
&=& \sum_{j\geq k-1}\P\left(S^{(b)}_{n-1}=j,S_n^{(b)}=j+1\right)\\
&=& \P\left(S^{(b)}_n\geq k\right).
\end{eqnarray*}
The last equality follows from the inclusion $\left\{S_n^{(b)}=j+1\right\}\subset\left\{S_{n-1}^{(b)}=j\right\}$ for any integer $j\geq 0$.
We then have shown that $\P\left(T_{n+1}\geq k\right)\geq \P\left(S_n^{(b)}\geq k\right)$ and (\ref{goodbound}) follows by induction on $n$.

Now from (\ref{goodbound}), we get that $\P\left(T_n=0\right)\leq \P\left(S_{n-1}^{(b)}=0\right)=b_{n-1}^{*}$ for $n\geq 1$ which completes the proof of the first point.

\item
This is immediate from the previous point since 
$$\P\left(\left(U_n^{x,y},\ldots,U_{n+k}^{x,y}\right)\neq \left(V_n^{x,y},\ldots,V_{n+k}^{x,y}\right)\right)\leq \sum_{j=n}^{n+k}\P\left(U_j^{x,y}\neq V_j^{x,y}\right).\square$$
\end{enumerate}

\begin{lem}\label{moment}
Let $(b_m)_{m\geq 0}$ be a sequence of non-negative real numbers decreasing to $0$ and such that $b_0<1$. Then for any nonnegative integer $k$, we have 
\begin{equation}\label{decayyy}
\sum_{m\geq 1}m^k b_m<\infty \Rightarrow \sum_{m\geq 1}m^k b_m^{*}<\infty.
\end{equation}
Moreover, if there exists $\rho\in (0,1)$ such that $b_m=O\left(\rho^m\right)$, then there exists $\overline{\rho}\in (0,1)$ s.t. $b_m^{*}=O\left(\overline{\rho}^m\right)$.
\end{lem}

\paragraph{Proof of Lemma \ref{moment}}
The case $k=0$ and the property of exponential decay has been treated in \citet{Bres}.
Let $(S_n)_{n\geq 0}$ be a Markov chain starting from $0$ and such that 
$$\P\left(S_{n+1}=i+1\vert S_n=i\right)=1-b_i=1-\P\left(S_{n+1}=0\vert S_n=i\right).$$
Let 
$$\tau=\inf\left\{n>0: S_n=0\right\}.$$
As explained in \citet{Bres}, equations $(5.7)$, $(5.11)$, $(A3)$ and $(A4)$, we have 
$$\P(\tau=1)=b_0,\quad \P(\tau=n)=b_{n-1}\prod_{m=0}^{n-2}(1-b_m),\quad n\geq 2$$
and the series $F(s)=\sum_{n\geq 1}\P\left(\tau=n\right)s^n$ and $G(s)=\sum_{n\geq 0}\P(S_n=0)$
satisfies $G(s)=(1-F(s))^{-1}$ for $s\in [0,1)$. 
We then have 
$$G'(s)=\frac{F'(s)}{(1-F(s))^2},\quad s\in [0,1).$$
From Beppo-Levi's theorem, we gave $F'(1^{-})=F'(1)$ which is finite using the assumption on the $b_m$'s and the bound $\P(\tau=n)\leq b_{n-1}$.

Moreover, using again Beppo Levi's theorem, we have $F(1^{-})=F(1)$. Moreover, as pointed out in \citet{Bres}, $(A.6)$, we have $F(1)<1$. We then conclude that $G'(1^{-})<\infty$. But once again from Beppo Levi's theorem, we have $G'(1^{-})=G'(1)=\sum_{n\geq 1}nb_n^{*}<\infty$.
This shows (\ref{decayyy}) for $k=1$. The case $k\geq 2$ is similar by taking the successive derivatives of $F$ and $G$. Details are omitted. $\square$

\begin{lem}\label{observation-driven}
Let $(X_t)_{t\in \Z}$ be a stationary process taking values in $\R^d$, $(y_t)_{t\in\Z}$ a sequence of point in $E$ and $\left\{G_{y,x}: (y,x)\in E\times \R^d\right\}$ a family of applications from $\R^k$ to $\R^k$ satisfying the three following assumptions.
\begin{enumerate}
\item
There exists $L\geq 1$ such that for $(y,y',x,x',z,z')\in E^2\times\R^{2d}\times \R^{2k}$, 
$$\left\vert G_{y,x}(z)-G_{y',x'}(z')\right\vert\leq L\left[\mathds{1}_{y\neq y'}+\vert x-x'\vert+\vert z-z'\vert\right].$$
\item
There exist a positive integer $r$ and a real number $\kappa\in (0,1)$ such that for all 
$(y_1,\ldots,y_r)\in E^r$, $(x_1,\ldots,x_r)\in \R^{dr}$ and $(s,s')\in \R^{2k}$,
$$\left\vert G_{y_1,x_1}\circ\cdots\circ G_{y_r,x_r}(s)-G_{y_1,x_1}\circ\cdots\circ G_{y_r,x_r}(s')\right\vert\leq \kappa\vert s-s'\vert.$$
\item
$\E\log^{+}\vert X_1\vert<\infty$.
\end{enumerate}
Then the following conclusions hold true.
\begin{itemize}
\item
Setting $\mathcal{D}=\left\{x\in (\R^d)^{\N}: \sum_{i=0}^{\infty}\kappa^{i/r}\vert x_i\vert<\infty\right\}$, we have $\P\left(X\in\mathcal{D}\right)=1$.
\item
For $t\in\Z$, $n\in \N^{*}$, $y\in E^{\Z}$ and $x\in\mathcal{D}$, set $\lambda_{n,t}^{(y,x)}=G_{y_{t-1},x_t}\circ\cdots\circ G_{y_{t-n-1},x_{t-n}}(0)$. Then for any $t\in\Z$, the sequence $\left(\lambda^{(y,x)}_{n,t}\right)_{n\geq 1}$ converges to an element of $\R^k$ denoted by $\lambda^{(y,x)}_t$. Moreover there exists $H:E^{\N}\times \mathcal{D}\rightarrow \R^k$ such that $\lambda_t=H\left(y_{t-1}^{-},x_t^{-}\right)$.
\item
Let $y\in E^{\Z}$ and $x\in (\R^k)^{\Z}$. If $\left(\overline{\lambda}_t\right)_{t\in\Z}$ is a sequence in $\R^k$, such that $\lim\inf_{t\rightarrow -\infty}\left\vert\overline{\lambda}_t\right\vert<\infty$ and $\overline{\lambda}_t=G_{y_{t-1},x_t}\left(\overline{\lambda}_{t-1}\right)$ for all $t\in\Z$, then $\overline{\lambda}_t=\lambda_t^{(y,x)}$ for all $t\in\Z$, where $\lambda_t^{(x,y)}$ is defined in the previous question. 
\item
Keeping the notations given in the previous point, we have for an integer $m\geq 1$
$$\sup_{x\in \mathcal{D}}\sup_{y_i=y'_i,-m+1\leq i\leq 0}\left\vert \lambda_0^{(y,x)}-\lambda_0^{(y',x)}\right\vert=O\left(\kappa^{m/r}\right).$$

\item
We have 
$$\left\vert H\left(y_{t-1}^{-},x_t^{-}\right)-H\left(y_{t-1}^{-},\overline{x}_t^{-}\right)\right\vert\leq \sum_{j=0}^{\infty}\kappa^j\sum_{i=1}^r L^i\left\vert x_{t-jr-i+1}-\overline{x}_{t-jr-i+1}
\right\vert.$$

\item
We assume that $\left((Y_t,X_t)\right)_{t\in\Z}$ is a stationary process taking values in $E\times \R^d$.
Then a stationary process $(\lambda_t)_{t\in\Z}$ satisfies the recursions $\lambda_t=G_{Y_{t-1},X_t}\left(\lambda_{t-1}\right)$ if and only if $\lambda_t=H\left(Y_{t-1}^{-},X_t^{-}\right)$, where $H$ is defined in the second point.

\end{itemize}
\end{lem}

\paragraph{Proof of Lemma \ref{observation-driven}}
We prove the result point by point.
\begin{itemize}
\item
The first point is a consequence of the following property. It $\beta\in (0,1)$, we have $\lim_{n\rightarrow -\infty}\beta^n X_n=0$ a.s. A proof of this fact can be found in \citet{DDM}, Lemma $34$.
\item
For $u\leq t$, we set $G_u^t=G_{y_{t-1},x_t}\circ\cdots\circ G_{y_u,x_u}$ and $G_t=G_t^t$. For $s,s'\in\R^k$, we have, using the two first assumptions, 
$$\left\vert G_{t-n}^t(0)-G_{t-n-1}^t(0)\right\vert\leq \kappa^{\frac{n+1}{r}-1}L^r\left\vert G_{t-n-1}^{t-n-1}(0)\right\vert.$$
If $(\overline{y},\overline{x})$ denotes a reference point in $E\times \R^d$, we also have
$$\left\vert G_{t-n-1}^{t-n-1}(0)-G_{\overline{y},\overline{x}}(0)\right\vert\leq L\left(1+\vert x_{t-n-1}-\overline{x}\vert\right).$$
Since $x_0^{-}\in\mathcal{D}$ and $\kappa\in (0,1)$, this clearly shows that $\sum_{n=0}^{\infty}\left\vert G_{t-n}^t(0)-G_{t-n-1}^t(0)\right\vert<\infty$. 
We then deduce that $\lim_{n\rightarrow \infty}G_{t-n}^t(0)$ exists. The existence of a measurable function $H$ follows, with $H\left(y_{t-1}^{-},x_t^{-}\right)=\lim_{n\rightarrow\infty}H_n\left(y_{t-1}^{-},x_t^{-}\right)$ and 
$H_n\left(y_{t-1}^{-},x_t^{-}\right)=G_{t-n}^t(0)$.

\item
Let $\left(\overline{\lambda}_t\right)_{t\in\Z}$ a sequence satisfying the proposed assumptions. Let $t\in\Z$. There exists a sequence $(n_i)_{i\in\N}$ of positive integers, such that $\lim_{i\rightarrow \infty}n_i=\infty$ and the sequence $\left(\overline{\lambda}_{t-n_i-1}\right)_{i\in \N}$ is bounded.
Writing $\overline{\lambda}_t=G_{t-n_i}^t\left(\overline{\lambda}_{t-n_i-1}\right)$, we have
$$\left\vert G_{t-n_i}^t(0)-\overline{\lambda}_t\right\vert L^r \kappa^{\frac{n_i+1}{r}-1}\left\vert \lambda_{t-n_i-1}\right\vert.$$
Letting $i\rightarrow \infty$, we deduce that $\lambda_t^{y,x}=\lambda_t$.

\item
We denote by $\overline{y}\in E^{\N}$ an arbitrary sequence.
We will first bound $\left\vert \lambda^{y,x}_t-\lambda^{\overline{y},x}_t\right\vert$. 
To this end we set, for $t\in\Z$, $\overline{G}_t=G_{\overline{y}_t,x_t}$. Let also $s\in\R^k$.
First, note that from our two first assumptions, we have
\begin{eqnarray*}
\left\vert G_{t-r+1}^t(s)-\overline{G}_{t-r+1}^t(s)\right\vert&\leq&  L\left\vert G_{t-r+1}^{t-1}(s)-\overline{G}_{t-r+1}^t(s)\right\vert+ \left\vert G_t\circ \overline{G}_{t-r+1}^{t-1}-\overline{G}_{t-r+1}^t(s)\right\vert\\
&\leq& L\left\vert G_{t-r+1}^{t-1}(s)-\overline{G}_{t-r+1}^t(s)\right\vert +L\\
\end{eqnarray*}
We then deduce that $\left\vert G_{t-r+1}^t(s)-\overline{G}_{t-r+1}^t(s)\right\vert\leq \sum_{i=1}^r L^i:=\widetilde{L}$ and the previous bound does not depend on $s$. 
Next, we get for a positive integer $n$,
$$\left\vert G_{t-nr+1}^t(0)-\overline{G}_{t-nr+1}^t(0)\right\vert \kappa \left\vert G_{t-nr+1}^{t-r}(0)-\overline{G}_{t-nr+1}^{t-r}(0)\right\vert + \widetilde{L}.$$ 
We then deduce that
$$\left\vert G_{t-nr+1}^t(0)-\overline{G}_{t-nr+1}^t(0)\right\vert \leq \frac{\widetilde{L}}{1-\kappa}.$$
Letting $n\rightarrow \infty$, we get $\left\vert \lambda_t^{(y,x)}-\overline{\lambda}_t^{(\overline{y},x)}\right\vert\leq \frac{\widetilde{L}}{1-\kappa}$.

Now, assume that the sequence $\overline{y}$ is such that $\overline{y}_i=y_i$, $0\leq i\leq m-1$. 
We have 
\begin{eqnarray*}
\left\vert \lambda_0^{(y,x)}-\lambda_0^{(\overline{y},x)}\right\vert &=& \left\vert G_{-m+1}^0\left(\lambda_{-m}^{(y,x)}\right)-G_{-m+1}^0\left(\lambda_{-m}^{(\overline{y},x)}\right)\right\vert\\&\leq& L^r\kappa^{\frac{m}{r}-1}\left\vert \lambda_{-m}^{(y,x)}-\lambda_{-m}^{(\overline{y},x)}\right\vert\\
&\leq&  L^r\kappa^{\frac{m}{r}-1}\cdot \frac{\widetilde{L}}{1-\kappa}.
\end{eqnarray*}
The result follows from the previous bounds.

\item 
The proof is very similar to that of the previous point and is then omitted.

\item
Assume first that $(\lambda_t)_{t\in\Z}$ is stationary and satisfies $\lambda_t=G_{Y_{t-1},X_t}\left(\lambda_{t-1}\right)$ a.s. In particular, we have $\lim\inf_{t\rightarrow -\infty}\left\vert \lambda_t\right\vert<\infty$ a.s. Hence, from the second point, we have $\lambda_t=\lambda^{(Y,X)}_t$ a.s.
On the other hand, suppose that $\lambda_t=H\left(Y_{t-1}^{-},X_t^{-}\right)$. We then have $\lambda_t=\lim_{n\rightarrow \infty}
\lambda_{n,t}^{(Y,X)}$ a.s. Since $G_{Y_{t-1},X_t}\left(\lambda_{n,t-1}^{(Y,X)}\right)=\lambda_{n+1,t}^{(Y,X)}$,
and the $G_{y,x}'$s are continuous, letting $n\rightarrow \infty$, we get $\lambda_t=G_{Y_{t-1},X_t}\left(\lambda_{t-1}\right)$ a.s. The equivalence between the two assertions then follows.$\square$

\end{itemize}

\bibliographystyle{plainnat}
\bibliography{biblogistic}

\begin{thebibliography}{44}
\providecommand{\natexlab}[1]{#1}
\providecommand{\url}[1]{\texttt{#1}}
\expandafter\ifx\csname urlstyle\endcsname\relax
  \providecommand{\doi}[1]{doi: #1}\else
  \providecommand{\doi}{doi: \begingroup \urlstyle{rm}\Url}\fi

\bibitem[Bressaud et~al.(1999{\natexlab{a}})Bressaud, Fern\'andez, and
  Galves]{Bres}
X.~Bressaud, R.~Fern\'andez, and A.~Galves.
\newblock Decay of correlations for non-{H}\"{o}lderian dynamics. a coupling
  approach.
\newblock \emph{Electron. J. Probab.}, 4:\penalty0 1--19, 1999{\natexlab{a}}.

\bibitem[Bressaud et~al.(1999{\natexlab{b}})Bressaud, Fern{\'a}ndez, and
  Galves]{Bres2}
Xavier Bressaud, Roberto Fern{\'a}ndez, and Antonio Galves.
\newblock Speed of d-convergence for markov approximations of chains with
  complete connections. a coupling approach.
\newblock \emph{Stochastic processes and their applications}, 83\penalty0
  (1):\penalty0 127--138, 1999{\natexlab{b}}.

\bibitem[Busch et~al.(2009)Busch, Ferrari, Flesia, Fraiman, Grynberg, and
  Leonardi]{prot}
Jorge~R Busch, Pablo~A Ferrari, Ana~Georgina Flesia, Ricardo Fraiman,
  Sebastian~P Grynberg, and Florencia Leonardi.
\newblock Testing statistical hypothesis on random trees and applications to
  the protein classification problem.
\newblock \emph{The Annals of Applied Statistics}, 3\penalty0 (2):\penalty0
  542--563, 2009.

\bibitem[Candelon et~al.(2013)Candelon, Dumitrescu, Hurlin, and
  Palm]{Candelson}
Bertrand Candelon, Elena-Ivona Dumitrescu, Christophe Hurlin, and Franz~C Palm.
\newblock Multivariate dynamic probit models: an application to financial
  crises mutation.
\newblock In \emph{VAR Models in Macroeconomics--New Developments and
  Applications: Essays in Honor of Christopher A. Sims}, pages 395--427.
  Emerald Group Publishing Limited, 2013.

\bibitem[Chamberlain(1982)]{Chamb}
Gary Chamberlain.
\newblock The general equivalence of granger and sims causality.
\newblock \emph{Econometrica: Journal of the Econometric Society}, pages
  569--581, 1982.

\bibitem[Cogburn(1984)]{Cog}
R.~Cogburn.
\newblock The ergodic theory of {M}arkov chains in randon environments.
\newblock \emph{Z. Wahrscheinlichkeitstheory verw. Gebiete}, 66:\penalty0
  109--128, 1984.

\bibitem[Comets et~al.(2002)Comets, Fern{\'a}ndez, Ferrari, et~al.]{comets}
Francis Comets, Roberto Fern{\'a}ndez, Pablo~A Ferrari, et~al.
\newblock Processes with long memory: regenerative construction and perfect
  simulation.
\newblock \emph{The Annals of Applied Probability}, 12\penalty0 (3):\penalty0
  921--943, 2002.

\bibitem[Cox et~al.(1981)Cox, Gudmundsson, Lindgren, Bondesson, Harsaae, Laake,
  Juselius, and Lauritzen]{Cox}
David~R Cox, Gudmundur Gudmundsson, Georg Lindgren, Lennart Bondesson, Erik
  Harsaae, Petter Laake, Katarina Juselius, and Steffen~L Lauritzen.
\newblock Statistical analysis of time series: Some recent developments [with
  discussion and reply].
\newblock \emph{Scandinavian Journal of Statistics}, pages 93--115, 1981.

\bibitem[de~Jong and Woutersen(2011)]{deJong2011}
Robert~M de~Jong and Tiemen Woutersen.
\newblock Dynamic time series binary choice.
\newblock \emph{Econometric Theory}, 27:\penalty0 673--702, 2011.

\bibitem[Dedecker and Prieur(2004)]{DedP1}
J{\'e}r{\^o}me Dedecker and Cl{\'e}mentine Prieur.
\newblock Coupling for $\tau$-dependent sequences and applications.
\newblock \emph{Journal of Theoretical Probability}, 17\penalty0 (4):\penalty0
  861--885, 2004.

\bibitem[Dedecker and Prieur(2005)]{DedP2}
J{\'e}r{\^o}me Dedecker and Cl{\'e}mentine Prieur.
\newblock New dependence coefficients. examples and applications to statistics.
\newblock \emph{Probability Theory and Related Fields}, 132\penalty0
  (2):\penalty0 203--236, 2005.

\bibitem[Dedecker et~al.(2007)Dedecker, Doukhan, Lang, Rafael, Louhichi, and
  Prieur]{DandCo}
J{\'e}r{\^o}me Dedecker, Paul Doukhan, Gabriel Lang, Le{\'o}n R~Jos{\'e}
  Rafael, Sana Louhichi, and Cl{\'e}mentine Prieur.
\newblock Weak dependence.
\newblock In \emph{Weak dependence: With examples and applications}, pages
  9--20. Springer, 2007.

\bibitem[Doeblin and Fortet(1937)]{DF}
W.~Doeblin and R.~Fortet.
\newblock Sur les cha\^{i}nes \`a liaisons compl\`etes.
\newblock \emph{Bull. Soc. Math. France}, 65:\penalty0 132--148, 1937.

\bibitem[Douc et~al.(2013)Douc, Doukhan, and Moulines]{DDM}
Randal Douc, Paul Doukhan, and Eric Moulines.
\newblock Ergodicity of observation-driven time series models and consistency
  of the maximum likelihood estimator.
\newblock \emph{Stochastic Processes and their Applications}, 123\penalty0
  (7):\penalty0 2620--2647, 2013.

\bibitem[Doukhan(1994)]{Doukhan(1994)}
P.~Doukhan.
\newblock \emph{Mixing: properties and examples}.
\newblock Number~85 in Lecture Notes in Statistics. Springer-Verlag, New York,
  1994.

\bibitem[Eichler et~al.()Eichler, Manner, and Turk]{Eichler}
Michael Eichler, Hans Manner, and Dennis Turk.
\newblock Dynamic copula based multivariate discrete choice models.
\newblock
  \emph{http://citeseerx.ist.psu.edu/viewdoc/download?doi=10.1.1.707.4242\&rep=rep1\&type=pdf}.

\bibitem[Fernandez and Galves(2002)]{GF}
Roseli Fernandez and A~Galves.
\newblock Markov approximations of chains of infinite order.
\newblock \emph{Bulletin of the Brazilian Mathematical Society}, 33\penalty0
  (3):\penalty0 295--306, 2002.

\bibitem[Fokianos and Truquet(2019)]{FT1}
K.~Fokianos and L.~Truquet.
\newblock On categorical time series with covariates.
\newblock \emph{Stochastic processes and their applications}, 129:\penalty0
  3446--3462, 2019.

\bibitem[Fokianos et~al.(2009)Fokianos, Rahbek, and Tjostheim]{Fok}
K.~Fokianos, A.~Rahbek, and D.~Tjostheim.
\newblock Poisson autoregression.
\newblock \emph{J. Amer. Statist. Assoc.}, 104:\penalty0 1430--1439, 2009.

\bibitem[Fokianos and Kedem(2003)]{Fokianos}
Konstantinos Fokianos and Benjamin Kedem.
\newblock Regression theory for categorical time series.
\newblock \emph{Statist. Sci.}, 18:\penalty0 357--376, 2003.
\newblock ISSN 0883-4237.
\newblock \doi{10.1214/ss/1076102425}.
\newblock URL \url{http://dx.doi.org/10.1214/ss/1076102425}.

\bibitem[Galves et~al.(2012)Galves, Galves, Garcia, Garcia, and
  Leonardi]{Galves}
Antonio Galves, Charlotte Galves, Jesus~E Garcia, Nancy~L Garcia, and Florencia
  Leonardi.
\newblock Context tree selection and linguistic rhythm retrieval from written
  texts.
\newblock \emph{The Annals of Applied Statistics}, 6\penalty0 (1):\penalty0
  186--209, 2012.

\bibitem[Guanche et~al.(2014)Guanche, M{\'\i}nguez, and M{\'e}ndez]{Guanche}
Yanira Guanche, Roberto M{\'\i}nguez, and Fernando~J M{\'e}ndez.
\newblock Autoregressive logistic regression applied to atmospheric circulation
  patterns.
\newblock \emph{Climate dynamics}, 42\penalty0 (1-2):\penalty0 537--552, 2014.

\bibitem[Hao et~al.(2016)Hao, Hao, Xia, Singh, Hong, Shen, and Ouyang]{Hao}
Zengchao Hao, Fanghua Hao, Youlong Xia, Vijay~P Singh, Yang Hong, Xinyi Shen,
  and Wei Ouyang.
\newblock A statistical method for categorical drought prediction based on
  nldas-2.
\newblock \emph{Journal of Applied Meteorology and Climatology}, 55\penalty0
  (4):\penalty0 1049--1061, 2016.

\bibitem[Harris(1955)]{Harris}
T.E. Harris.
\newblock On chains of infinite order.
\newblock \emph{Pacific J. Math.}, 5:\penalty0 707--724, 1955.

\bibitem[Iosifescu and Grigorescu(1990)]{Ios}
M.~Iosifescu and S.~Grigorescu.
\newblock \emph{Dependence with Complete Connections and its Applications}.
\newblock Cambridge University Press, 1990.

\bibitem[Kaufmann(1987)]{Kaufmann}
H.~Kaufmann.
\newblock Regression models for nonstationary categorical time series:
  Asymptotic estimation theory.
\newblock \emph{Annals of Statistics}, 15:\penalty0 79--98, 1987.

\bibitem[Kauppi and Saikkonen(2008)]{kauppi}
Heikki Kauppi and Pentti Saikkonen.
\newblock Predicting us recessions with dynamic binary response models.
\newblock \emph{The Review of Economics and Statistics}, 90\penalty0
  (4):\penalty0 777--791, 2008.

\bibitem[Kifer(1996)]{Kif}
Yuri Kifer.
\newblock Perron-frobenius theorem, large deviations, and random perturbations
  in random environments.
\newblock \emph{Mathematische Zeitschrift}, 222\penalty0 (4):\penalty0
  677--698, 1996.

\bibitem[Lindvall(2002)]{Lind}
Torgny Lindvall.
\newblock \emph{Lectures on the coupling method}.
\newblock Courier Corporation, 2002.

\bibitem[L{\"u}tkepohl(2005)]{Lut}
H.~L{\"u}tkepohl.
\newblock \emph{New Introduction to Multiple Time Series Analysis}.
\newblock Springer, Berlin, 1st edition, 2005.

\bibitem[Merlev{\`e}de et~al.(2011)Merlev{\`e}de, Peligrad, and Rio]{Mer}
Florence Merlev{\`e}de, Magda Peligrad, and Emmanuel Rio.
\newblock A bernstein type inequality and moderate deviations for weakly
  dependent sequences.
\newblock \emph{Probability Theory and Related Fields}, 151\penalty0
  (3-4):\penalty0 435--474, 2011.

\bibitem[Mitrophanov(2005)]{MIT}
A~Yu Mitrophanov.
\newblock Sensitivity and convergence of uniformly ergodic markov chains.
\newblock \emph{Journal of Applied Probability}, 42\penalty0 (4):\penalty0
  1003--1014, 2005.

\bibitem[Neumann(2011)]{Neu}
M.~Neumann.
\newblock Absolute regularity and ergodicity of poisson count processes.
\newblock \emph{Bernoulli}, 17:\penalty0 1268--1284, 2011.

\bibitem[Nyberg(2014)]{Nyberg}
Henri Nyberg.
\newblock A bivariate autoregressive probit model: Business cycle linkages and
  transmission of recession probabilities.
\newblock \emph{Macroeconomic Dynamics}, 18\penalty0 (4):\penalty0 838--862,
  2014.

\bibitem[Park et~al.(2017)Park, Simar, and Zelenyuk]{Park}
Byeong~U Park, L{\'e}opold Simar, and Valentin Zelenyuk.
\newblock Nonparametric estimation of dynamic discrete choice models for time
  series data.
\newblock \emph{Computational Statistics \& Data Analysis}, 108:\penalty0
  97--120, 2017.

\bibitem[Rao(2009)]{PR}
BLS~Prakasa Rao.
\newblock Conditional independence, conditional mixing and conditional
  association.
\newblock \emph{Annals of the Institute of Statistical Mathematics},
  61\penalty0 (2):\penalty0 441--460, 2009.

\bibitem[Rissanen(1983)]{Ris}
Jorma Rissanen.
\newblock A universal data compression system.
\newblock \emph{IEEE Transactions on information theory}, 29\penalty0
  (5):\penalty0 656--664, 1983.

\bibitem[Russell and Engle(2005)]{russell}
Jeffrey~R Russell and Robert~F Engle.
\newblock A discrete-state continuous-time model of financial transactions
  prices and times: The autoregressive conditional multinomial--autoregressive
  conditional duration model.
\newblock \emph{Journal of Business \& Economic Statistics}, 23\penalty0
  (2):\penalty0 166--180, 2005.

\bibitem[Rydberg and Shephard(2003)]{rydberg}
Tina~Hviid Rydberg and Neil Shephard.
\newblock Dynamics of trade-by-trade price movements: decomposition and models.
\newblock \emph{Journal of Financial Econometrics}, 1\penalty0 (1):\penalty0
  2--25, 2003.

\bibitem[Sims(1972)]{Sims}
Christopher~A Sims.
\newblock Money, income, and causality.
\newblock \emph{The American economic review}, 62\penalty0 (4):\penalty0
  540--552, 1972.

\bibitem[Villani(2009)]{Villani}
C.~Villani.
\newblock \emph{Optimal Transport. Old and New}.
\newblock Springer, 2009.

\bibitem[Woodard et~al.(2011)Woodard, Matteson, Henderson, et~al.]{EJS}
Dawn~B Woodard, David~S Matteson, Shane~G Henderson, et~al.
\newblock Stationarity of generalized autoregressive moving average models.
\newblock \emph{Electronic Journal of Statistics}, 5:\penalty0 800--828, 2011.

\bibitem[Wu(2005)]{Wu}
Wei~Biao Wu.
\newblock Nonlinear system theory: Another look at dependence.
\newblock \emph{Proceedings of the National Academy of Sciences}, 102\penalty0
  (40):\penalty0 14150--14154, 2005.

\bibitem[Yuan and Lei(2013)]{CC}
DeMei Yuan and Lan Lei.
\newblock Some conditional results for conditionally strong mixing sequences of
  random variables.
\newblock \emph{Science China Mathematics}, 56\penalty0 (4):\penalty0 845--859,
  2013.

\end{thebibliography}

\end{document}